\documentclass[11pt,reqno]{amsart}
\usepackage[T1]{fontenc}
\usepackage[utf8]{inputenc}
\usepackage{lmodern}
\usepackage[a4paper,margin=27mm]{geometry}
\usepackage{mathtools,amssymb,mathrsfs,xcolor}
\usepackage{float}
\usepackage{bm}
\usepackage{microtype}
\usepackage[colorlinks,
breaklinks,
linkcolor=red!80!black,
anchorcolor=green,
citecolor=blue!60!black, 
]{hyperref}

\usepackage{enumitem}
\setlist[enumerate]{itemsep=0.2em,topsep=0.3em,leftmargin=2em}

\newtheorem{theorem}{Theorem}[section]
\newtheorem{proposition}[theorem]{Proposition}
\newtheorem{lemma}[theorem]{Lemma}

\theoremstyle{remark}
\newtheorem{remark}{Remark}[section]
\theoremstyle{definition}
\newtheorem{defn}{Definition}[section]
\theoremstyle{plain}
\newtheorem*{question}{Question}

\newcommand{\R}{\mathbb R}
\newcommand{\Z}{\mathbb Z}
\newcommand{\T}{\mathbb T}
\newcommand{\cA}{\mathcal A}
\newcommand{\cAt}{\widetilde{\mathcal A}}

\newcommand{\tcI}{\widetilde{\mathcal I}}
\newcommand{\cI}{\mathcal I}
\newcommand{\cP}{\mathcal P}
\newcommand{\cU}{\mathcal U}
\newcommand{\Lip}{\operatorname{Lip}}

\newcommand{\supp}{\operatorname{supp}}
\newcommand{\meas}{\mathcal L^1}

\newcommand{\vep}{\varepsilon}
\newcommand{\hL}{\mathbb{L}}

\newcommand{\nAct}{A_{\mathbb{L}}}
\newcommand{\bgamma}{\bm{\gamma}}
\newcommand{\dbgamma}{\dot{\bm{\gamma}}}
\newcommand{\bT}{\mathbf{T}}
\newcommand{\kt}{h_t^\mathbb{L}}
\newcommand{\kT}{h_T^\mathbb{L}}
\newcommand{\kinfty}{h^\mathbb{L}}

\def\tilde{\widetilde}
\def\hat{\widehat}

\title[]{Hyperbolic periodic orbits and optimal $1/t$ convergence of the Lax--Oleinik semigroup}
\author{Qinbo Chen}
\address[Qinbo Chen]{School of Mathematics, Nanjing University, Nanjing 210093, China}
\email{qinbochen@nju.edu.cn}

\author{Zhi-Xiang Zhu}
\address[Zhi-Xiang Zhu]{School of Mathematics, Nanjing University, Nanjing 210093, China}
\email{zhuzhixiang0928@gmail.com}
\keywords{Lax--Oleinik semigroup, Hamilton--Jacobi equations, Aubry set, hyperbolicity}

\begin{document}

\begin{abstract}
Let $M$ be a closed manifold and $H:T^*M\to\mathbb{R}$ a $C^k$ ($k\geq2$) Tonelli Hamiltonian whose Aubry set consists of finitely many hyperbolic orbits, each of which is either a periodic orbit or a fixed point. We construct a $C^k$ critical subsolution that is quadratically strict outside the projected Aubry set. We then use this subsolution to establish an $O(1/t)$ rate of convergence for the associated Lax--Oleinik semigroup and the minimal action functions. This rate is optimal whenever at least one of these hyperbolic orbits is periodic. 
\end{abstract}

\maketitle
%\tableofcontents

\section{Introduction and Main Results}
\label{sec_introduction}

Let $M$ be a closed, connected Riemannian manifold of dimension $n$, and let $T^*M$ and $TM$ denote its cotangent and tangent bundles, respectively. Throughout this paper, we consider a {\it Tonelli Hamiltonian} $H:T^*M\to\R$ of class $C^k$ ($k\geq 2$). That is, $H\in C^k(T^*M)$ satisfies
\begin{enumerate}
	\item  [1)]{\bf Strict convexity.} The fiberwise Hessian $\partial_{pp}^2 H(x,p)$ is positive definite for every $(x,p)\in T^*M$.
	
	\item  [2)]{\bf Superlinearity.} For every $x\in M$, $H(x,p)/\|p\|_x\to+\infty$ as $\|p\|_x\to\infty$. 
\end{enumerate}
The Legendre--Fenchel transform of $H$ gives rise to an associated Lagrangian $L:TM\to\R$. It is well known (see, e.g., \cite{Fathi_book}) that $L$ is also of class $C^k$, strictly convex and superlinear in each fiber $T_xM$. Thus, $L$ is a {\it Tonelli Lagrangian}. 

For $t>0$ and $x,y\in M$, let $h_t(x,y)$ denote the minimal action in time $t$, defined by
\begin{equation*}
h_t(x,y)
:=\inf_{\substack{  \gamma(0)=x, \gamma(t)=y}}\,\int_0^t L(\gamma(s),\dot\gamma(s)) \,ds,
\end{equation*}
where the infimum is taken over continuous piecewise $C^1$ curves. For $u\in C(M)$ and $t> 0$, we define $T_t^-u:M\to\R$ by
\[
T_t^-u(x):=\inf_{y\in M}\left\{u(y)+h_t(y,x)\right\}.
\]
Setting $T_0^-u=u$, we obtain the family of nonlinear operators $\{T_t^-\}_{t\geq 0}$, known as the {\it Lax--Oleinik semigroup}. The study of its properties is important in several areas, including PDEs, dynamical systems, optimization and control theory.

The evolution $U(x,t):=T_t^-u(x)$ generated by this semigroup is the unique viscosity solution of the Cauchy problem for the Hamilton--Jacobi equation
\begin{equation}\label{eq_eHJ}
 U_t+H\bigl(x, d_xU\bigr)=0 \quad\text{in }M\times(0,\infty), \qquad U(x,0)=u(x).
\end{equation}
The large-time behavior of solutions to \eqref{eq_eHJ} has been extensively studied. Early contributions were made by Namah and Roquejoffre \cite{Namah_Roquejoffre1999} using PDE methods and by Fathi \cite{Fathi98_laxoleinik} using a dynamical systems approach, while among the many subsequent works we mention \cite{Barles_Souganidis2000, Roquejoffre2001, Davini_Siconolfi2006, Ishii2008_Asymptotic}. In particular, Fathi \cite{Fathi98_laxoleinik} proved that for every $u\in C(M)$, the limit
\[
u_\infty(x):=\lim_{t\to +\infty}\big(T_t^-u(x)+c_Ht\big)
\]
exists uniformly on $M$, where $c_H$ denotes {\it Ma\~n\'e's critical value}, also known as the {\it ergodic constant}. The limiting function $u_\infty$ is a viscosity solution of the stationary Hamilton--Jacobi equation
\begin{equation}\label{eq_sHJ}
 H(x,d_xv)=c_H\quad\text{in }M.
\end{equation}
In weak KAM theory \cite{Fathi_book}, $u_\infty$ is also a backward weak KAM solution. Further geometric perspectives on the Lax--Oleinik semigroup were provided by Arnaud \cite{Arnaud_05LaxOleinik, Arnaud_2011Pseudographs} (see also \cite{Bernard_C11, Bernard08_pseudographs}).

In this paper, we study the rate of convergence of $T_t^-u+c_Ht$ to $u_\infty$ as $t\to \infty$. Of particular interest is how the dynamical properties of the underlying system influence this convergence rate. Previous works (see, for instance, \cite{Iturriaga_Morgado09, Wang_Yan11Scichina, Wang_Yan12CMP, Wang_Yan2012nonli, Morgado15_hyperbolicity, Wangkaizhi18, Iturriaga_Khanin_Zhang20}) suggest that this rate is closely related to the dynamical structure of the Aubry set (see Definition \ref{defn_Aubry}). According to a famous conjecture of Ma\~n\'e, the Aubry set of a generic Tonelli Hamiltonian consists of either a hyperbolic fixed point or a hyperbolic periodic orbit \cite{Mane1996} (see also \cite{Contreras_Figalli_Rifford_15, Contreras24} for related discussion on surfaces). Motivated by these considerations, we therefore investigate the rate of convergence of $T_t^-u+c_Ht$ to $u_\infty$ under a hyperbolicity assumption on the Aubry set.

When the Aubry set consists of hyperbolic fixed points, Iturriaga and S\'anchez-Morgado \cite{Iturriaga_Morgado09} established the following exponential convergence estimate: there exists $\mu>0$ such that, for every $u\in C(M)$, there exists $K>0$ satisfying
\begin{equation}\label{IM_exp}
	\|T_t^-u+c_Ht-u_\infty\|_{C(M)} \le K e^{-\mu t},\qquad t\geq 0.
\end{equation}
In contrast, when the Aubry set consists of hyperbolic periodic orbits, the convergence need not be exponential.  Indeed, Wang and Yan \cite[Example 4.1]{Wang_Yan2012nonli} constructed a particularly instructive example of an autonomous Tonelli Lagrangian on $M=\T^2$:
\begin{equation}\label{wangyanexample}
L(x_1,x_2,v_1,v_2)=\frac12v_1^2+1-\cos(2\pi x_1)+\frac12(v_2-a)^2, \qquad a\neq0.
\end{equation}
In this example, the Aubry set consists of a single periodic orbit $\cAt=\{(0,y,0,a):y\in\T\}$, which is hyperbolic on its energy level. They showed that there exist an initial datum $u\in C(\T^2)$, a point $q_0\in\T^2$, a constant $C>0$, and a sequence of times $t_n\to+\infty$ such that
\[
|T_{t_n}^-u(q_0)+c_H t_n-u_\infty(q_0)|
\geq Ct_n^{-1}.
\]
Wang and Yan's example shows that, even when the Aubry set reduces to a single hyperbolic periodic orbit, a convergence rate faster than order $1/t$ cannot be expected in general. 

This naturally raises the following question: 

\begin{question}
For every autonomous Tonelli system whose Aubry set consists of finitely many hyperbolic periodic orbits, does the estimate $\|T_t^-u+c_Ht-u_\infty\|_{C(M)}=O(1/t)$ as $t\to +\infty$ hold for all $u\in C(M)$?  If so, is this rate optimal for each such system?
\end{question}

We answer both questions affirmatively in Theorems \ref{Thm_semigroup} and \ref{Thm_optimality}.

\begin{theorem}\label{Thm_semigroup}
Let $H:T^*M\to\R$ be a $C^k$ Tonelli Hamiltonian ($k\geq2$), and assume that its Aubry set consists of finitely many hyperbolic orbits, each of which is either a periodic orbit or a fixed point. Then there exist constants $C>0$ and $t_*>0$ such that, for any $u\in C(M)$,
\begin{equation*}
\left\|T_t^-u+c_Ht-u_\infty\right\|_{C(M)}\le\frac{C}{t},
\qquad \text{for all } t\ge t_*.
\end{equation*}
\end{theorem}

\begin{remark}
Throughout this paper, a periodic orbit is always understood to be nontrivial (i.e., not a fixed point). For an autonomous Hamiltonian or Euler--Lagrange flow, the hyperbolicity of a periodic orbit must be understood with respect to the flow restricted to its energy level. See Subsection \ref{subsection_hyperbolicity}. 
\end{remark}

The $O(1/t)$ rate in Theorem \ref{Thm_semigroup} can be improved to an exponential rate when the Aubry set consists entirely of hyperbolic fixed points (see \eqref{IM_exp} above). In contrast, we show that the $O(1/t)$ rate is \textbf{optimal} whenever the Aubry set contains at least one periodic orbit.

\begin{theorem}\label{Thm_optimality}
Under the hypotheses of Theorem \ref{Thm_semigroup}, suppose further that the Aubry set contains at least one hyperbolic periodic orbit. Then there exist an initial datum $u\in C^k(M)$, constants $0<c<C$, and $t_*>0$ such that
\[
\frac{c}{t}\leq \left\|T_t^-u+c_Ht-u_\infty\right\|_{C(M)}\leq \frac{C}{t},
\qquad  \text{for all }t\ge t_*.
\]	
\end{theorem}

Theorem \ref{Thm_optimality} shows that the $1/t$ lower-bound behavior observed in the example \eqref{wangyanexample} from \cite{Wang_Yan2012nonli} also occurs for any autonomous Tonelli system whose Aubry set consists of finitely many hyperbolic periodic orbits.

It is worth noting that the lower bound in Theorem \ref{Thm_optimality} cannot hold for every initial datum. Indeed, any backward weak KAM solution $u$ satisfies $T_t^-u+c_Ht=u=u_\infty$ for all $t\geq0$. 

Theorem \ref{Thm_semigroup} is, in fact, a direct consequence of the following convergence estimate for the minimal action functions. To state it, recall that the Peierls barrier of $L$ is defined by
\begin{equation*}
h(x,y):=\liminf_{t\to\infty}\big( h_t(x,y)+c_Ht\big).
\end{equation*}
See Section \ref{sec_prelimi} for further details. We establish the following result.

\begin{theorem}\label{Thm_kernel}
Under the hypotheses of Theorem \ref{Thm_semigroup}, there exist $C>0$ and $t_*>0$ such that
\begin{equation}\label{h_upperbd}
\left\|h_t+c_Ht-h\right\|_{C(M\times M)}\leq \,\frac{C}{t}, \qquad \text{for all } t\ge t_*.
\end{equation}
Suppose further that the Aubry set contains at least one hyperbolic periodic orbit. Then there exists $c>0$ such that
\begin{equation}\label{h_lowerbd}
 \left\|h_t+c_Ht-h\right\|_{C(M\times M)}\geq \,\frac{c}{t}, \qquad \text{for all } t\ge t_*.
\end{equation}
\end{theorem}

A key ingredient in the proofs of Theorems \ref{Thm_semigroup}--\ref{Thm_kernel} is the construction of a $C^k$ critical subsolution that is quadratically strict outside the projected Aubry set $\cA\subset M$. See Theorem \ref{Thm_subsol_quad} below. 

The Hamilton--Jacobi equation \eqref{eq_sHJ} need not admit a $C^1$ solution. Fathi and Siconolfi \cite{Fathi_Siconolfi2004} established the existence of $C^1$ critical subsolutions that are strict outside $\cA$. More precisely, there exists $u_0\in C^1(M)$ such that  $H(x,d_xu_0)=c_H$ for $x\in\cA$, and
\[
H(x,d_xu_0)<c_H\qquad\text{for }x\in M\setminus\cA.
\]
Bernard \cite{Bernard_C11} subsequently improved the regularity to $C^{1,1}$. Furthermore, Bernard proved in \cite{Bernard_smoothsol07} that such a subsolution $u_0$ can be chosen to be as smooth as the Tonelli Hamiltonian when the Aubry set  consists of hyperbolic periodic orbits and hyperbolic fixed points.

We prove a quantitative refinement of Bernard's result by establishing a quadratic strictness estimate in terms of the distance to the projected Aubry set $\cA$.

\begin{theorem}[Smooth critical subsolution with quadratic strictness]
\label{Thm_subsol_quad}
Under the hypotheses of Theorem \ref{Thm_semigroup}, there exist $f\in C^k(M)$ and $\kappa>0$ such that $H(x,d_xf)=c_H$ for $x\in\cA$, and 
\[
H(x,d_xf)\leq c_H-\kappa\,d(x,\cA)^2,\quad \text{for }x\in M.
\]
Here, $d(x,\cA):=\inf_{y\in\cA}d(x,y)$ denotes the Riemannian distance from $x$ to $\cA$.
\end{theorem}

Although our focus is on autonomous systems, we briefly recall, for comparison, some related results for nonautonomous Hamiltonians $H(x,p,t)$. Fathi and Mather \cite{Fathi_Mather} showed that the Lax--Oleinik semigroup need not converge in the time-periodic setting. Wang and Yan \cite{Wang_Yan12CMP} introduced new L--O type operators and later proved their exponential convergence when the Aubry set is a single hyperbolic periodic orbit \cite{Wang_Yan2012nonli}. S\'anchez-Morgado \cite{Morgado15_hyperbolicity} proved exponential convergence to time-periodic solutions on the torus when the Aubry set consists of hyperbolic periodic orbits. A related result for closed manifolds was obtained in \cite{Wangkaizhi18}.	 The difference between the $O(1/t)$ convergence rate in the autonomous setting and the exponential convergence rate in the time-periodic setting is related to the role of neutral directions. In the autonomous case, a hyperbolic periodic orbit necessarily has a neutral flow direction in $T^*M$ (see \eqref{periodicsplitting}). In the time-periodic case, however, the flow direction along a hyperbolic periodic orbit lies in the extended phase space $T^*M\times\mathbb S^1$ and is transverse to each fixed-time section, so it does not induce a neutral direction for the return map on $T^*M$.

\noindent\textbf{Organization of the paper.} Section \ref{sec_prelimi} recalls the necessary background from weak KAM theory and hyperbolic dynamics. In Section \ref{sec_quadraticsol}, we prove Theorem \ref{Thm_subsol_quad} and then apply it in Section \ref{sec_connections} to establish some quantitative action estimates for connecting curves near the Aubry set. These constructions are used in Section \ref{sec_convergencerate} to establish the upper-bound estimate \eqref{h_upperbd} in Theorem \ref{Thm_kernel} and Theorem \ref{Thm_semigroup}. Finally, Section \ref{sec_optimal} proves the optimality result in Theorem \ref{Thm_optimality} and the estimate \eqref{h_lowerbd} in Theorem \ref{Thm_kernel}.

\section{Preliminaries}
\label{sec_prelimi}
Throughout this section, let $M$ be a closed Riemannian manifold of dimension $n$, and let $H: T^*M\to \R$ be an autonomous $C^k$ ($k\ge 2$) Tonelli Hamiltonian. The Lagrangian $L: TM\to\R$ associated with $H$ is 
\[
L(x,v)=\max_{p\in T_x^*M}\, \{\langle p, v \rangle_x-H(x,p)\},\qquad \forall~(x,v)\in TM. 
\]
Then $L$ is a $C^k$ Tonelli Lagrangian. 

\subsection{Weak KAM theory and Aubry sets}
We recall some useful results from Fathi's weak KAM theory for Tonelli systems  \cite{Fathi_book}. We also refer the reader, for instance, to \cite{Contreras_Iturriaga_book, Sorrentino_book, ChengXuebook, Zavidovique_book2025} for further discussions of weak KAM theory and its applications.

Under the above Tonelli assumptions, one can define a $C^{k-1}$ diffeomorphism between $TM$ and $T^*M$, called the {\it Legendre transform}:
\begin{equation*}
\mathscr L ~:~ TM \longrightarrow T^*M,\qquad
(x,v)  \longmapsto   \left(x,\frac{\partial L}{\partial v}(x,v)\right).
\end{equation*}
Let $\Phi_t^L$ and $\Phi_t^H$ denote the Euler--Lagrange flow and the Hamiltonian flow, respectively. The two $C^{k-1}$ flows are equivalent from a dynamical systems point of view, in the sense that they are conjugate via the Legendre transform, namely, 
\begin{equation}\label{flowconjugacy}
\Phi_t^L=\mathscr{L}^{-1}\circ\Phi_t^H\circ\mathscr L.
\end{equation}

The (stationary) Hamilton--Jacobi equation associated with $H$ is 
\begin{equation}\label{HJSansu}
H(x,d_xu)=c \quad\text{in}~ M.
\end{equation} 
It is well known (see, e.g., \cite{Fathi2012, Contreras_Iturriaga_book}) that there exists a unique constant $c=c_H$, called the {\it critical value}, for which equation \eqref{HJSansu} admits solutions in the viscosity sense.

The Lagrangian action $A_L(\xi)$ of a continuous piecewise $C^1$ curve $\xi:[a,b]\to M$ is
\begin{equation*}
	A_L(\xi):=\int_a^b L(\xi(s),\dot\xi(s))\,ds.
\end{equation*}
A continuous piecewise $C^1$ curve $\gamma:[a,b]\to M$ is called a {\it minimizing curve} (for $L$) if $A_L(\gamma)\le A_L(\xi)$ holds for any continuous piecewise $C^1$ curve $\xi: [a,b]\to M$ with $\xi(a)=\gamma(a)$, $\xi(b)=\gamma(b)$. By Tonelli's theorem (see, e.g., \cite[Section 3]{Fathi_book}), we have

\begin{proposition}\label{minimizer_regularity}
	If $\gamma:[a,b]\to M$ is a minimizing curve, then it is of class $C^k$.
\end{proposition}

Following Fathi \cite{Fathi_book}, we recall the notions of domination and calibration. 
\begin{defn}
Let $u:M\to\R$ and $c\in\R$.
\begin{enumerate}
\item  We say that $u$ is {\it dominated by $L+c$}, and write $u\prec L+c$, if for  every continuous piecewise $C^1$ curve $\xi:[t_1,t_2]\to M$,
\[
u(\xi(t_2))-u(\xi(t_1))\le A_L(\xi)+c(t_2-t_1).
\]

\item  Suppose that $u\prec L+c$. A continuous piecewise $C^1$ curve $\gamma: I\to M$ is said to be {\it $(u, L, c)$-calibrated} if
\[
u(\gamma(t_2))-u(\gamma(t_1))=A_L\left (\gamma|_{[t_1,t_2]}\right)+c(t_2-t_1),\quad \text{for any}~t_1<t_2\in I.
\]
\end{enumerate}
\end{defn}

Every calibrated curve is necessarily a minimizing curve for $L$ and hence is of class $C^k$. In addition, $u\prec L+c$ if and only if $u$ is a viscosity subsolution of $H(x,d_xu)=c$. A fundamental feature of dominated functions, as observed by Fathi \cite{Fathi_book, Fathi14ICM}, is that they naturally give rise to invariant sets of the Euler--Lagrange flow $\Phi_t^L$. More precisely, given a viscosity subsolution $u$ of 
\begin{equation}\label{eqsHJc_H}
	H(x, d_xu)=c_H \quad\text{in}~ M,
\end{equation}
equivalently, $u\prec L+c_H$, we define the subset $\tcI(u)$ of $TM$ by
\begin{equation*}
	\tcI(u)
	:=
	\left\{
	 (x,v)\in TM ~:~t\mapsto\pi\circ \Phi_t^L(x,v) \text{~is $(u,L,c_H)$-calibrated on~}\R 
	 \right\},
\end{equation*}
where $\pi: TM\to M$ is the canonical projection. 

\begin{proposition}\label{prop_u_calib}
$\tcI(u)$ is a nonempty compact $\Phi_t^L$-invariant subset of $TM$. For every $(x_0,v_0)\in \tcI(u)$,  the
function $u$ is differentiable at $x_0$ and
\begin{equation*}
	d_{x_0}u=\partial_v L(x_0,v_0),\qquad    H\left(x_0,\partial_v L(x_0,v_0) \right)=c_H.
\end{equation*}
\end{proposition}

Via the Legendre transform, we define the corresponding subset of $T^*M$ by 
\[\tcI^*(u)=\mathscr L\left(\tcI(u)\right).\] 
By Proposition \ref{prop_u_calib}, we have $H\big|_{\tcI^*(u)}=c_H$.

\begin{defn}\label{defn_Aubry}
The {\it Aubry set} $\cAt$ of $L$ is defined by 
\begin{equation*}
	\cAt:=\bigcap_{u\in \mathcal{SS}^-}\,\tcI(u) 
\end{equation*}
where $\mathcal{SS}^-$ is the set of all viscosity subsolutions of \eqref{eqsHJc_H}. Via the Legendre transform, we define
\[\cAt^*:=\mathscr L\left(\cAt\right)=\bigcap_{u\in \mathcal{SS}^-}\,\tcI^*(u), \]
which is called the Aubry set of $H$.
\end{defn}

We write 
\[\cI(u):=\pi\left( \tcI(u)\right),\qquad \cA:=\pi\left(\cAt\right)\]
for the projections of $\tcI(u)$ and $\cAt$ onto $M$, respectively.  We call $\cA$ the {\it projected Aubry set}.

\begin{proposition}[Mather's graph property \cite{Mather1993, Fathi_book}]\label{prop_Matherproperty}
	The Aubry set $\cAt$  (resp. $\cAt^*$) is a nonempty compact subset of $TM$  (resp. $T^*M$), and is invariant under the flow $\Phi_t^L$ (resp. $\Phi_t^H$). Moreover, for every  viscosity subsolution $u$ of \eqref{eqsHJc_H}, 
	\[ 
	\tcI^*(u)\,\subset \,\operatorname{Graph}(du):=\{(x, d_xu)\,:\, x\in M \text{~such that $d_x u$ exists} \}.
	\]
	In particular, $\cAt$ (resp. $\cAt^*$) is a Lipschitz graph over the projected Aubry set $\cA$, namely, $\pi|_{\cAt}:\cAt\to\cA$ is a bi-Lipschitz homeomorphism.
\end{proposition}

For $t>0$, let $h_t$ denote the minimal action in time $t$, defined by
\[
h_t(x,y):=\inf_{\xi\in \Gamma^t_{x,y}}A_L(\xi), 
\]
where $\Gamma^t_{x,y}$ is the family of all continuous piecewise $C^1$ curves $\xi:[0,t]\to M$ satisfying $\xi(0)=x$ and $\xi(t)=y$. By Tonelli's theorem, the infimum is attained. 

\begin{proposition}\label{Prop_Fleming}
For each $t_0>0$, the family $\{h_t+c_Ht\}_{t\ge t_0}$ is uniformly bounded on $M\times M$.
\end{proposition}
See \cite[Section 5]{Fathi_book} for its proof.  We next recall the definition of the Peierls barrier \cite{Mather1993}.

\begin{defn}
The {\it Peierls barrier} of $L$ is the function $h:M\times M\to \R$ defined by
\begin{equation*}
	h(x,y):=\liminf_{t\to+\infty}\, \big(h_t(x,y)+c_Ht\big).
\end{equation*}	
\end{defn}
In the autonomous case, it was proved in \cite{Fathi98_laxoleinik} that $h_t(x,y)+c_Ht$ actually converges uniformly to $h(x,y)$ on $M\times M$. The function $h$ is Lipschitz continuous. Moreover, the projected Aubry set satisfies $\cA=\{x\in M: h(x,x)=0\}$.

For $t>0$, the Lax--Oleinik operator 
$T_t^-: C(M)\to C(M)$
can be expressed using $h_t$ by
\begin{equation*}
	T_t^-u(x)=\inf_{y\in M} \{u(y)+h_t(y,x)\}.
\end{equation*}
The semigroup $u\mapsto T_t^-u+c_Ht$, $t\geq0$, admits a fixed point \cite{Fathi97}. The following theorem describes its long-time behavior.

\begin{theorem}[\cite{Fathi98_laxoleinik}]\label{uinftyexpression}
For every $u\in C(M)$, the uniform limit $u_\infty$ $:=$ $\lim_{t\to+\infty}(T_t^-u+c_Ht)$ exists. It satisfies $T_t^-u_\infty+c_Ht=u_\infty$ for all $t\geq 0$. Moreover,
\[
u_\infty(x)=\inf_{y\in M}\bigl\{u(y)+h(y,x)\bigr\},
\qquad x\in M.
\]
\end{theorem}

Following Fathi, we recall the definition of weak KAM solutions. 

\begin{defn}
Let $u\in C(M)$ be dominated by $L+c_H$.
\begin{enumerate}
\item $u$ is a {\it backward weak KAM solution} of $L$ if for each $x\in M$ there exists a $(u,L, c_H)$-calibrated curve $\gamma_x:(-\infty,0]\to M$ such that $\gamma_x(0)=x$;

\item $u$ is a {\it forward weak KAM solution} of $L$ if for each $x\in M$ there exists a $(u,L,c_H)$-calibrated curve $\gamma_x:[0,+\infty)\to M$ such that $\gamma_x(0)=x$.
\end{enumerate}	
\end{defn}

Weak KAM solutions are Lipschitz continuous on $M$ but may have singularities. For the role of Lax--Oleinik operators in singularity propagation, see, for instance, \cite{Cannarsa_Cheng17CVPDE, Cannarsa_Chen_ChengJDE, Cannarsa_Cheng_Fathi21, CCHW24}. However, weak KAM solutions are differentiable on the projected Aubry set $\cA$.

\begin{proposition}[\cite{Fathi_book, Fathi_Siconolfi2004}]
\label{prop_commondifferential}
Let $u^-$ and $u^+$ be arbitrary backward and forward weak KAM solutions, respectively. Then they are both differentiable at every $x\in\cA$, and
\[
(x,d_xu^-)=(x,d_xu^+)\in \cAt^*, \qquad x\in\cA.
\]
In particular, $H(x, d_xu^-)=H(x, d_xu^+)=c_H$ for all $x\in\cA$. 

Moreover, there exists a backward weak KAM solution $u$ such that $\tcI^*(u)=\cAt^*$.
\end{proposition}

We finally describe the relationship between the Lax--Oleinik semigroup $T_t^-$, backward weak KAM solutions, and viscosity solutions. We refer to \cite{Fathi_book, Fathi2012} for its proof.
\begin{theorem}\label{calibrationThm}
For $w\in C(M)$, the following statements are equivalent:
(i) $w$ is a backward weak KAM solution; 
(ii) $T_t^-w+c_Ht=w$ for all $t\geq 0$;
(iii) $w$ is a viscosity solution of \eqref{eqsHJc_H}.
\end{theorem}

\subsection{Hyperbolic dynamics}\label{subsection_hyperbolicity}
We recall some basic facts from hyperbolic dynamics. For more details, we refer to the standard references \cite{Hirsch_Pugh_Shub_book, Katokbook}. Recall that $\Phi_t^H$ denotes the $C^{k-1}$ Hamiltonian flow on $T^*M$, generated by the Hamiltonian vector field $X_H$.

\begin{defn}\label{defn_hyperb_orbits}
 A fixed point $z_0$ of $\Phi_t^H$ is said to be {\it hyperbolic} if, for every $t\neq 0$, $D\Phi_t^H(z_0):T_{z_0}(T^*M)\to T_{z_0}(T^*M)$ has no eigenvalues on the unit circle. A (nontrivial) periodic orbit $\cP$ of the flow $\Phi_t^H$, lying on the energy level $\Sigma_e:=H^{-1}(e)$, is said to be {\it hyperbolic on its energy level} if the linearized Poincar\'e return map, taken on a local section transverse to the flow within $\Sigma_e$, has no eigenvalues on the unit circle.
\end{defn}

We remark that, by \eqref{flowconjugacy}, a fixed point or a periodic orbit of the Euler--Lagrange flow is hyperbolic if and only if its image under the Legendre transform $\mathscr L$ is hyperbolic for the Hamiltonian flow $\Phi_t^H$.

We now recall some standard properties of the local stable and unstable manifolds associated with hyperbolic fixed points and periodic orbits.

\smallskip
\noindent{\bf The fixed-point case.}
Let $z_0$ be a hyperbolic fixed point of $\Phi_t^H$. By hyperbolicity and the symplectic structure, the tangent space $T_{z_0}(T^*M)$ admits a $D\Phi_t^H$-invariant splitting
\[
T_{z_0}(T^*M)=E^s(z_0)\oplus E^u(z_0),
\qquad
\dim E^s(z_0)=\dim E^u(z_0)=n,
\]
and there exist constants $C>0$ and $\mu>0$ such that, for all $t\geq 0$,
\[
\|D\Phi_t^H(z_0)v\| \leq Ce^{-\mu t}\|v\|,\quad v\in E^s(z_0);\qquad
\|D\Phi_{-t}^H(z_0)v\| \leq Ce^{-\mu t}\|v\|,\quad v\in E^u(z_0).
\]

The stable manifold theorem yields $C^{k-1}$ embedded local stable and unstable manifolds $W^s_{\mathrm{loc}}(z_0)$ and $W^u_{\mathrm{loc}}(z_0)$. For sufficiently small $\delta>0$, they can be described as
\begin{equation*}
\begin{aligned}
W^s_{\mathrm{loc}}(z_0)
:=~&
\left\{
y\in T^*M:
d\left(\Phi_t^H(y),z_0\right)<\delta,\quad
\text{for all }t\geq 0
\right\},
\\
W^u_{\mathrm{loc}}(z_0)
:=~&
\left\{
y\in T^*M:
d\left(\Phi_t^H(y),z_0\right)<\delta,\quad
\text{for all }t\leq 0
\right\}.
\end{aligned}
\end{equation*}
Here, $d(\cdot,\cdot)$ denotes the distance induced by the Riemannian metric on $T^*M$. In fact, every point on these manifolds converges exponentially to $z_0$ in the corresponding time direction. Moreover, $T_{z_0}W^s_{\mathrm{loc}}(z_0)=E^s(z_0)$ and $T_{z_0}W^u_{\mathrm{loc}}(z_0)=E^u(z_0)$. In particular,
\begin{equation}\label{transverse_fixed}
T_{z_0}W^s_{\mathrm{loc}}(z_0)
\oplus T_{z_0}W^u_{\mathrm{loc}}(z_0)
=T_{z_0}(T^*M),
\qquad
T_{z_0}W^s_{\mathrm{loc}}(z_0)
\cap T_{z_0}W^u_{\mathrm{loc}}(z_0)
=\{0\}.
\end{equation}
Thus, although both local invariant manifolds lie in the energy level $H^{-1}(e)$ with $e=H(z_0)$, their transversality is understood in the full phase space $T^*M$.

\smallskip
\noindent{\bf The periodic-orbit case.} 
Let $\cP\subset\Sigma_e$ be a periodic orbit that is hyperbolic on its energy level $\Sigma_e:=H^{-1}(e)$. Let $T>0$ be its least period. Fix $\tilde{x}\in\cP$. Then $\cP=\{\Phi_t^H(\tilde{x}):t\in[0,T]\}$. Since $X_H(z)\neq 0$ for all $z\in\cP$, we have $d_zH\neq 0$. Hence, $\Sigma_e$ is a $C^k$ hypersurface near $\cP$ with 
\[
T_z\Sigma_e=\ker d_zH,\qquad \forall~z\in\cP.
\]

By hyperbolicity, there exist continuous $D\Phi_t^H$-invariant subbundles $E^s=\bigsqcup_{z\in\cP}E^s(z)$ and $E^u=\bigsqcup_{z\in\cP}E^u(z)$ along $\cP$, such that, for every $z\in\cP$,
\begin{equation}\label{periodicsplitting}
	T_z\Sigma_e=E^s(z)\oplus \R X_H(z)\oplus E^u(z),\qquad \dim E^s(z)=\dim E^u(z)=n-1.
\end{equation}
Moreover, $D\Phi_t^H$ contracts $E^s$ exponentially in forward time and $E^u$ exponentially in backward time, uniformly along $\cP$.

Applying the stable manifold theorem to the flow on $\Sigma_e$, we obtain $C^{k-1}$ embedded local stable and unstable manifolds $W^s_{\mathrm{loc}}(\cP)$ and $W^u_{\mathrm{loc}}(\cP)$, both of dimension $n$.  For sufficiently small $\delta>0$, they can be described as
\begin{equation*}
\begin{aligned}
W^s_{\mathrm{loc}}(\cP)
:=~&
\left\{
y\in\Sigma_e:
d\left(\Phi_t^H(y),\cP\right)<\delta,\quad
\text{for all }t\geq 0
\right\},
\\
W^u_{\mathrm{loc}}(\cP)
:=~&
\left\{
y\in\Sigma_e:
d\left(\Phi_t^H(y),\cP\right)<\delta,\quad
\text{for all }t\leq 0
\right\}.
\end{aligned}
\end{equation*}
In fact, points in $W^s_{\mathrm{loc}}(\cP)$ converge exponentially to $\cP$ in forward time, while points in $W^u_{\mathrm{loc}}(\cP)$ converge exponentially to $\cP$ in backward time. Moreover, for every $z\in\cP$,
\[
T_zW^s_{\mathrm{loc}}(\cP)
=E^s(z)\oplus\R X_H(z),
\qquad
T_zW^u_{\mathrm{loc}}(\cP)
=E^u(z)\oplus\R X_H(z).
\]
Consequently, the two manifolds intersect transversely along $\cP$ within the energy level $\Sigma_e$, in the sense that, for every $z\in\cP$,
\begin{equation}\label{transverse}
T_zW^s_{\mathrm{loc}}(\cP)
+
T_zW^u_{\mathrm{loc}}(\cP)
=
T_z\Sigma_e,
\qquad
T_zW^s_{\mathrm{loc}}(\cP)
\cap
T_zW^u_{\mathrm{loc}}(\cP)
=
\R X_H(z).
\end{equation}

Thus, by \eqref{transverse_fixed} and \eqref{transverse}, in both the fixed-point case and periodic-orbit case, we have
\begin{equation}\label{eq_transverse}
T_zW^s_{\mathrm{loc}}(\cP)
\cap T_zW^u_{\mathrm{loc}}(\cP)
=T_z\cP,
\qquad z\in\cP,
\end{equation}
where $\cP=\{z_0\}$ and $T_{z_0}\cP=\{0\}$ in the fixed-point case.

When $\cP$ is contained in the Aubry set $\cAt^*$ of $H$, the local stable and unstable manifolds admit the following graph representations. See, e.g., \cite{Contreras_Iturriaga_conjugate99, Bernard_smoothsol07}.

\begin{proposition}\label{prop_graphforWsu}
Let $\cP\subset\cAt^*$ be either a hyperbolic periodic orbit or a hyperbolic fixed point, where $\cP=\{z_0\}$ in the fixed-point case. Then there exist small neighborhoods $U^s$ and $U^u$ of $\pi(\cP)$ in $M$, together with $C^{k-1}$ closed $1$-forms $\omega^s$ on $U^s$ and $\omega^u$ on $U^u$, such that
\[
W^s_{\mathrm{loc}}(\cP)\cap\pi^{-1}(U^s)
=\{(x,\omega^s_x):x\in U^s\},
\qquad
W^u_{\mathrm{loc}}(\cP)\cap\pi^{-1}(U^u)
=\{(x,\omega^u_x):x\in U^u\}.
\]
Here, $\pi:T^*M\to M$ denotes the canonical projection. 
\end{proposition}

By Proposition \ref{prop_commondifferential}, there always exists a backward weak KAM solution $u$ such that $\tcI^*(u)=\cAt^*$. The following localization result follows from Bernard \cite{Bernard_smoothsol07}.

\begin{theorem}[{\cite[Section 6]{Bernard_smoothsol07}}] 
\label{thm_localization_calibrated_orbits}
Under the assumptions of Theorem \ref{Thm_semigroup}, let $u\in C(M)$ be a backward weak KAM solution satisfying $\tcI^*(u)=\cAt^*$, and let $\cP\subset\cAt^*$ be either a hyperbolic periodic orbit or a hyperbolic fixed point. 

Then there exist a small compact neighborhood $V$ of $\cP$ in $\Sigma_{c_H}=H^{-1}(c_H)$, with $V\cap\cAt^*=\cP$, and a small neighborhood $U$ of $\pi(\cP)$ in $M$, such that the following properties hold:
\begin{enumerate}
\item [\rm(i)] If $z\in V$ and $\Phi_t^H(z)\in V$ for all $t\leq 0$, then
$
z\in W^u_{\mathrm{loc}}(\cP).
$

\item[\rm(ii)] For any $(u, L, c_H)$-calibrated curve $\gamma:(-\infty,0]\to  M$ with $\gamma(0)\in U$,  its Hamiltonian lift $\tilde\gamma(t):= \mathscr L\big(\gamma(t),\dot\gamma(t)\big)$ remains in $V$ for all $t\leq 0$.

\item [\rm(iii)] Denote 
$
\Gamma_u:=\left\{
(x,p)\in T^*M:\,
t\mapsto\pi\circ\Phi_t^H(x,p)
\text{~is $(u, L, c_H)$-calibrated on~} (-\infty,0]
\right\}.
$
Then, 
\begin{equation}\label{Wuinclusion}
	\Gamma_u\cap\pi^{-1}(U)
\subset
W^u_{\mathrm{loc}}(\cP).
\end{equation}
Consequently, $u$ is of class $C^k$ on $U$, and 
\[
\operatorname{Graph}\left(du|_{U}\right)
=
W^u_{\mathrm{loc}}(\cP)
\cap\pi^{-1}(U)
\]
where $\operatorname{Graph}\left(du|_{U}\right)
=\{(x, d_xu): x\in U\}$.
\end{enumerate}
\end{theorem}

If $\cP$ is a periodic orbit, assertion (i) follows from the unstable manifold theorem applied to the flow on the energy level $\Sigma_{c_H}$ near $\cP$. If $\cP=\{z_0\}$ is a hyperbolic fixed point, apply the unstable manifold theorem in $T^*M$ and choose a sufficiently small ambient neighborhood $V_0$ of $z_0$. Setting $V=V_0\cap\Sigma_{c_H}$ gives (i) in this case as well.

It is known (see, e.g., \cite{Fathi_book, Fathi_Siconolfi2004}) that the Hamiltonian lift of every $(u,L,c_H)$-calibrated curve lies on the energy level $\Sigma_{c_H}$. Assertion (ii) then follows from \cite[Lemma 8]{Bernard_smoothsol07}. 
	
Assertion (iii) follows from \cite[Proposition 7]{Bernard_smoothsol07} and its proof. For the reader's convenience, we include the proof here. For each $z=(x,p)\in\Gamma_u\cap\pi^{-1}(U)$, by the definition of $\Gamma_u$ and assertion (ii), we have $\Phi_t^H(z)\in V$ for all $t\leq 0$. Hence, by assertion (i), $z\in W^u_{\mathrm{loc}}(\cP)$, which proves \eqref{Wuinclusion}. 

Shrinking $U$ if necessary, Proposition \ref{prop_graphforWsu} yields a $C^{k-1}$ closed $1$-form $\omega^u$ on $U$ such that 
\begin{equation}\label{eqwugraph}
W^u_{\mathrm{loc}}(\cP)\cap\pi^{-1}(U)
=
\{(x,\omega_x^u):x\in U\}.
\end{equation}
Since $u$ is a backward weak KAM solution, $\Gamma_u\cap T_x^*M\neq \emptyset$ for every $x\in M$. Combining this with \eqref{Wuinclusion} and the graph property \eqref{eqwugraph}, we obtain
\[
\Gamma_u\cap T_x^*M
=
\{(x,\omega_x^u)\},\quad \text{for each}~x\in U.
\]
By the standard differentiability criterion for backward weak KAM solutions, the set $\Gamma_u\cap T_x^*M$ is a singleton if and only if $u$ is
differentiable at $x$, in which case $\Gamma_u\cap T_x^*M=\{(x,d_xu)\}$. Consequently, $d_xu=\omega_x^u$ for all $x\in U$. Since $\omega^u$ is $C^{k-1}$,  $u\in C^k(U)$. Moreover,
\begin{equation*}
\operatorname{Graph}\left(du|_U\right)
=
\{(x,\omega_x^u):x\in U\}
=
W^u_{\mathrm{loc}}(\cP)\cap\pi^{-1}(U).
\end{equation*}

\section{Hyperbolicity and Smooth Critical Subsolutions}
\label{sec_quadraticsol}
\paragraph{Notation.}
In the following sections, we adopt the notation below:

\begin{itemize}
\item Let $\cAt^*\subset T^*M$ be the Aubry set of $H$. Under the hypothesis of Theorem \ref{Thm_semigroup}, we write
\begin{equation*}\label{hyp_HP}\tag{HP}
\cAt^*=\cAt_1^*\cup\cdots\cup\cAt_m^*,
\end{equation*}
where $m\geq1$ and the sets $\cAt_1^*,\ldots,\cAt_m^*$ are distinct connected components of $\cAt^*$. Each $\cAt_i^*$ is a hyperbolic periodic orbit or fixed point of the Hamiltonian flow.

\item The projected Aubry set $\cA\subset M$ is then given by
\[\cA=\bigcup_{i=1}^m\cA_i,\quad\text{with}\quad \cA_i:=\pi(\cAt^*_i),\]
where $\pi: T^*M\to M$ is the canonical projection. 
\end{itemize}

The main goal of this section is to prove Theorem \ref{Thm_quadratic}, from which Theorem \ref{Thm_subsol_quad} follows. 

Bernard \cite{Bernard_smoothsol07} established the following existence result for smooth critical subsolutions of \eqref{eqsHJc_H}.
\begin{theorem}[Bernard]\label{thm_smoothCksol}
Let $H$ be a $C^k$ Tonelli Hamiltonian satisfying the hyperbolicity condition \eqref{hyp_HP}.  Then there exists a $C^k$ critical subsolution $u_0: M\to\R$ of \eqref{eqsHJc_H} such that 
\[
H(x, d_xu_0)=c_H\quad \text{for}~ x\in\cA, 
\quad\text{and}\quad  
H(x,d_xu_0)<c_H
\quad \text{for}~ x\in M\setminus\cA.
\]
\end{theorem}

For our purposes, we establish a quantitative refinement of Theorem \ref{thm_smoothCksol}, with a quadratic estimate in terms of the distance to the projected Aubry set. 
\begin{theorem}\label{Thm_quadratic}
Under the assumptions of Theorem \ref{thm_smoothCksol}, there exists $f\in C^k(M)$ such that
\[
H(x,d_xf)
\leq
c_H-\kappa\,d(x,\cA)^2, \quad x\in M;\qquad\quad
H(x,d_xf)
=
c_H,\quad x\in\cA,
\]
where $\kappa>0$ is a constant. In particular, $\tcI^*(f)=\cAt^*$.
\end{theorem}

\subsection{Local regularity and differential separation}

\begin{lemma}
\label{lem_vpm}
There exist a backward weak KAM solution $v^-$ and a forward weak KAM solution $v^+$ on $M$, together with small open neighborhoods $U_i$ of $\cA_i$, $1\leq i\leq m$, such that, setting
\[
\cU:=\bigcup_{i=1}^m U_i,
\]
both $v^-$ and $v^+$ are $C^k$ on $\cU$, and, for each
$i=1,\ldots,m$,
\begin{equation}\label{vgraphproper}
\operatorname{Graph}\bigl(dv^-|_{U_i}\bigr)
=
W^u_{\mathrm{loc}}(\cAt_i^*)\cap\pi^{-1}(U_i),
\qquad
\operatorname{Graph}\bigl(dv^+|_{U_i}\bigr)
=
W^s_{\mathrm{loc}}(\cAt_i^*)\cap\pi^{-1}(U_i).
\end{equation}
\end{lemma}
See Figure \ref{fig_local_weak_kam_graphs} for a simple illustration of this local graph structure in the case $\dim M=2$.

\begin{proof}
We split the proof into two steps.

\smallskip
\noindent {\it Step 1}. We first construct the backward weak KAM solution. By Proposition \ref{prop_commondifferential}, we can choose a backward weak KAM solution $v^-:M\to\R$ satisfying $\tcI^*(v^-)=\cAt^*$.

Applying Theorem \ref{thm_localization_calibrated_orbits} to $v^-$ and each component $\cAt_i^*$, whether a hyperbolic fixed point or a hyperbolic periodic orbit, we obtain a small open neighborhood $O_i$ of $\cA_i$ such that
\[
v^-\in C^k(O_i),
\qquad
\operatorname{Graph}\bigl(dv^-|_{O_i}\bigr)
=
W^u_{\mathrm{loc}}(\cAt_i^*)\cap\pi^{-1}(O_i).
\]

\smallskip
\noindent {\it Step 2}. We next construct the forward weak KAM solution by time reversal, similarly to Fathi \cite[Section 4.9]{Fathi_book}. Define
\[
\check H(x,p):=H(x,-p),
\qquad
\check L(x,v):=L(x,-v),
\]
and let $\mathcal R:T^*M\to T^*M$, $\mathcal R(x,p)=(x,-p)$. Then $\check H$ is a $C^k$ Tonelli Hamiltonian with $c_{\check H}=c_H$, and its Hamiltonian flow satisfies
\begin{equation}\label{twoflows}
\Phi_t^{\check H}
=
\mathcal R\circ\Phi_{-t}^H\circ\mathcal R,
\qquad t\in\R.
\end{equation}
Moreover, the Aubry set $\check{\cA}^*$ of  $\check H$ satisfies 
\[
\check{\cA}^*
=
\bigcup_{i=1}^m \check{\cA}_i^*,
\qquad
\check{\cA}_i^*:=\mathcal R(\cAt_i^*),
\]
and each $\check{\cA}_i^*$ is a hyperbolic fixed point or a hyperbolic periodic orbit of $\check H$. By \eqref{twoflows}, the corresponding local unstable manifold $W^u_{\mathrm{loc}}(\check{\cA}_i^*)$ of $\check{\cA}_i^*$ satisfies 
\begin{equation}\label{cWutoWs}
\mathcal R\left(
W^u_{\mathrm{loc}}(\check{\cA}_i^*)
\right)
=
W^s_{\mathrm{loc}}(\cAt_i^*).
\end{equation}
Applying Step 1 to $\check H$, we obtain a backward weak KAM solution $\check v^-:M\to\R$ and open neighborhoods $\check O_i$ of $\cA_i$ such that
\begin{equation}\label{checkvminus}
\check v^-\in C^k(\check O_i),
\qquad
\operatorname{Graph}\bigl(d\check v^-|_{\check O_i}\bigr)
=
W^u_{\mathrm{loc}}(\check{\cA}_i^*)
\cap\pi^{-1}(\check O_i).
\end{equation}

Set
\[
v^+(x):=-\check v^-(x)\quad \text{for all}~x\in M.
\]
It is known in weak KAM theory (see \cite[Section 4.9]{Fathi_book}) that  $v^+$ is a forward weak KAM solution for $H$. Clearly, $v^+\in C^k(\check O_i)$, and $
\operatorname{Graph}\bigl(dv^+|_{\check O_i}\bigr)
=
\mathcal R\left(
\operatorname{Graph}\bigl(d\check v^-|_{\check O_i}\bigr)
\right)$. Combining this with \eqref{cWutoWs} and \eqref{checkvminus} yields 
\[
\operatorname{Graph}\bigl(dv^+|_{\check O_i}\bigr)
=
W^s_{\mathrm{loc}}(\cAt_i^*)
\cap\pi^{-1}(\check O_i).
\]

Finally, taking $U_i:=O_i\cap\check O_i$ and $\cU:=\bigcup_{i=1}^m U_i$, we finish the proof.
\end{proof}

\begingroup
\setlength{\intextsep}{1pt}
\setlength{\abovecaptionskip}{3pt}
\setlength{\belowcaptionskip}{0pt}
\begin{figure}[H]
  \centering
  \includegraphics[scale=1.2]{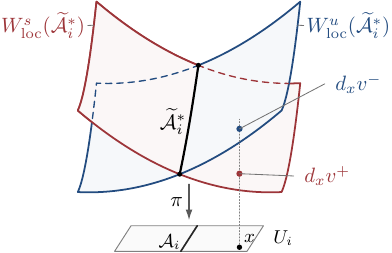}
  \caption{Local unstable and stable manifolds over $U_i$ as graphs of $dv^-$ and $dv^+$.}
  \label{fig_local_weak_kam_graphs}
\end{figure}
\endgroup

The transversality of the local stable and unstable manifolds yields a quantitative lower bound for $|dv^- - dv^+|$ near the projected Aubry set $\cA$.

\begin{lemma}
\label{lem:differential_gap}
Let $v^-$, $v^+$ be as given in Lemma \ref{lem_vpm}. After shrinking the neighborhoods $U_i$ if necessary, there exists a constant $\beta>0$ such that
\begin{equation}\label{eq_diffgap}
\left|d_xv^--d_xv^+\right|
\geq\beta\,d(x,\cA),
\qquad x\in\cU:=\bigcup_{i=1}^m U_i.
\end{equation}
Here $|\cdot|$ denotes the norm on $T_x^*M$ induced by the Riemannian metric, and $d(x,\cA):=\inf_{y\in\cA}d(x,y)$ denotes the distance from $x$ to $\cA$.
\end{lemma}

\begin{proof}
In view of Lemma \ref{lem_vpm}, we set
\[
\bm{h}:=v^--v^+.
\]
Then $\bm{h}\in C^k(\cU)$ where $\cU=\bigcup_{i=1}^m U_i$. By Proposition \ref{prop_commondifferential}, $d_xv^-=d_xv^+$ for all $x\in \cA$ and hence $d\bm{h}=0$ on $\cA$. Thus the Hessian $D^2\bm{h}(x)$ is intrinsically defined for every $x\in \cA$. We identify it with the linear map $D^2\bm{h}(x):T_xM\to T_x^*M$, $\xi\mapsto D^2\bm{h}(x)(\xi,\cdot)$. This implies
\[
\ker D^2\bm{h}(x)
=
\left\{
\xi\in T_xM:\,
D^2\bm{h}(x)(\xi,\eta)=0
\text{ for all }\eta\in T_xM
\right\}.
\]
Since $\cA_i$ is an embedded $C^k$ submanifold of $M$, we claim that
\begin{equation}\label{eq_hessiankernel}
\ker D^2\bm{h}(x)=T_x\cA_i, \quad\text{for every}~ x\in \cA_i.
\end{equation}
Here, if $\cA_i=\{a_i\}$ is a fixed point, then $T_{a_i}\cA_i=\{0\}$ and hence $D^2\bm{h}(a_i)$ is invertible. 

To prove \eqref{eq_hessiankernel}, we fix $x\in \cA_i$ and set
\[
z=(x,d_xv^-)=(x,d_xv^+)\in\cAt_i^*.
\]
At $z$, the graph identities \eqref{vgraphproper} give, in local cotangent coordinates,
\[
T_zW^u_{\mathrm{loc}}(\cAt_i^*)
=\left\{\big(\xi,D^2v^-(x)\xi\big): \xi\in T_xM \right\},
\qquad
T_zW^s_{\mathrm{loc}}(\cAt_i^*)
=\left\{\big(\xi,D^2v^+(x)\xi\big): \xi\in T_xM \right\}.
\]
Thus a vector $(\xi,D^2v^-(x)\xi)$ belongs to the intersection 
$
T_zW^u_{\mathrm{loc}}(\cAt_i^*)
\cap
T_zW^s_{\mathrm{loc}}(\cAt_i^*)
$
if and only if 
$\big(\xi,D^2v^-(x)\xi\big)=\big(\xi,D^2v^+(x)\xi\big)$. This is equivalent to
\[
D^2\bm{h}(x)\xi=0.
\]
Since the differential of the bundle projection
$\pi:T^*M\to M$ satisfies $D\pi_z(\xi,\eta)=\xi$, we get
\begin{equation}\label{eq_intersectionkernel}
D\pi_z\left(\,
T_zW^u_{\mathrm{loc}}(\cAt_i^*)
\cap
T_zW^s_{\mathrm{loc}}(\cAt_i^*)\,
\right)
=
\ker D^2\bm{h}(x).
\end{equation}
On the other hand, in view of the property \eqref{eq_transverse}, we obtain
\[
T_zW^u_{\mathrm{loc}}(\cAt_i^*)\cap T_zW^s_{\mathrm{loc}}(\cAt_i^*)=T_z\cAt_i^*.
\]
Combining this with \eqref{eq_intersectionkernel} yields $\ker D^2\bm{h}(x)=D\pi_z\bigl(T_z\cAt_i^*\bigr)=T_x\cA_i$, which verifies \eqref{eq_hessiankernel}.

We are now ready to prove \eqref{eq_diffgap}. It suffices to show that, for each
$i=1,\ldots,m$, there exists $\beta_i>0$ such that
\begin{equation}\label{dxh_ai}
|d_x\bm{h}|\geq \beta_i\,d(x,\cA_i),
\qquad x\in U_i,
\end{equation}
after shrinking $U_i$ if necessary.

\smallskip

{\it (i) Fixed-point case.} We first consider the case where $\cA_i=\{a_i\}$ corresponds to a fixed point. By shrinking $U_i$ if necessary, we choose a smooth coordinate chart $\chi:B_r^n\to U_i$ with $\chi(0)=a_i$, where $B_r^n=\{y\in\R^n:|y|<r\}$, and set
\[
g(y):=\bm{h}(\chi(y)).
\]
Then $Dg(0)=0$, and $B:=D^2g(0)$ is invertible by \eqref{eq_hessiankernel}. By Taylor's expansion, $Dg(y)=By+o(|y|)$. Since $B$ is invertible, shrinking $r>0$ if necessary, we obtain
\begin{equation}\label{eq_asiuq}
|Dg(y)|\geq c_0|y|, \qquad y\in B_r^n,
\end{equation}
for some $c_0>0$. We next translate \eqref{eq_asiuq} into an intrinsic estimate on $M$. By the local equivalence of norms and distances, there exists $C\geq 1$ such that
\[
|Dg(y)|\leq C|d_{\chi(y)}\bm{h}|,
\qquad
d(\chi(y),a_i)\leq C|y|.
\]
Consequently,
\[
|d_x\bm{h}|\geq c_0 C^{-2}\,d(x, a_i),
\qquad x\in U_i.
\]
Thus \eqref{dxh_ai} follows with $\beta_i=c_0C^{-2}$.

\smallskip

{\it (ii) Periodic-orbit case.} We next consider the case where $\cA_i$ corresponds to a periodic orbit. Observe that if $n=1$ where $n=\dim M$, then $\cA_i=M$, so the estimate is trivial. We may therefore assume $n\geq2$.

Fix $x_0\in \cA_i$. Since $\cA_i$ is an embedded $C^k$ circle, there exist $r>0$, a small neighborhood $O_{x_0}$ of $x_0$, and a $C^k$ coordinate chart
$
\chi:(-r,r)\times B_r^{n-1}\longrightarrow O_{x_0}
$
adapted to $\cA_i$, such that $\chi(0,0)=x_0$ and $\cA_i\cap O_{x_0}=\left\{\chi(s,0):\,|s|<r\right\}$. Set
\[
g(s,y):=\bm{h}\big(\chi(s,y)\big).
\]
$Dg(s,0)=0$ since $d_x\bm{h}=0$ for all $x\in \cA_i$. Moreover, \eqref{eq_hessiankernel} gives
\begin{equation}\label{kerd2f00}
	\ker D^2g(0,0)=\R(1,0).
\end{equation}
This, together with the symmetry of the Hessian, implies that
\[
D^2g(0,0)
=
\begin{pmatrix}
0 & 0\\
0 & B
\end{pmatrix}
\quad 
\text{where}~ 
B:=D^2_{yy}g(0,0).
\]
The $(n-1)\times (n-1)$ matrix $B$ is invertible by \eqref{kerd2f00}. Since $\partial_y g(s,0)=0$, 
\[
\partial_y g(s,y)
=By+\int_0^1\bigl(D^2_{yy} g(s,ty)-B\bigr)y\,dt
=By+o(|y|),\qquad\text{as}~(s,y)\to(0,0),
\]
where the last equality follows from the continuity of $D^2_{yy}g$. Since $B$ is invertible, shrinking $r>0$ if necessary, there exists $c_1>0$ such that
\begin{equation}\label{eq_gapcoordinate_estimate}
|\partial_y g(s,y)|\geq c_1|y|,
\qquad (s,y)\in(-r,r)\times B_r^{n-1}.
\end{equation}
By the local equivalence of norms and distances, there exists $K\geq1$ such that, for  $x=\chi(s,y)$ 
\[
|\partial_y g(s,y)|
\leq
K|d_x\bm{h}|,\qquad 
d(x,\cA_i)\leq K|y|.
\]
Combining these inequalities with \eqref{eq_gapcoordinate_estimate} yields
\begin{equation}\label{temp_est}
|d_x \bm{h}|
\geq
\frac{c_1}{K^2}\, d(x, \cA_i),\quad \text{for every } x\in O_{x_0}.
\end{equation}

Since $\cA_i$ is compact, finitely many coordinate neighborhoods of the above form cover $\cA_i$, each satisfying an estimate analogous to \eqref{temp_est}. Therefore, taking the minimum of the corresponding positive constants, we obtain \eqref{dxh_ai} after shrinking $U_i$ if necessary.

Finally, set $\beta:=\min_{1\leq i\leq m}\beta_i$. For any $x\in\cU=\cup_{i=1}^m U_i$, choose $i$ such that $x\in U_i$. Then
\[
|d_x\bm{h}| \geq \beta_i\,d(x,\cA_i) \geq \beta\,d(x,\cA).
\]
This proves \eqref{eq_diffgap}.
\end{proof}

\subsection{Proof of Theorem \ref{Thm_quadratic}}

\begin{proof}
Let $v^-:M\to\R$, $v^+:M\to \R$ and $\cU$ be as in Lemmas \ref{lem_vpm} and \ref{lem:differential_gap}. Hence, 
\begin{equation}\label{eq_dxvpmbeta}
|d_xv^--d_xv^+|
\geq
\beta\,d(x,\cA),
\qquad x\in\cU.
\end{equation}
Define $u_1: M\to\R$ by
\begin{equation}\label{defu1}
u_1:=\frac{1}{2}(v^-+v^+).
\end{equation}
Obviously, $u_1 \in C^k(\cU)$. Both $v^-$ and $v^+$ are viscosity subsolutions of $H(x,d_xu)=c_H$. The $p$-convexity of $H$ thus implies that $u_1$ is also a viscosity subsolution. 

\medskip

\noindent{\bf Step 1: local quadratic estimate.} We first claim that there is a constant $a>0$ such that
\begin{equation}\label{eq_quadraticgap}
H(x,d_xu_1)\leq c_H-a\,d(x,\cA)^2,
\qquad x\in\cU.
\end{equation}

Indeed, the graph properties \eqref{vgraphproper} imply that
\begin{equation}\label{Hxvpm_cH}
H(x,d_xv^-)=H(x,d_xv^+)=c_H, \qquad x\in\cU.
\end{equation}
Since $v^-$ and $v^+$ are Lipschitz on $M$, there exists $K>0$ such that $\Lip(v^\pm)\leq K$. Hence, $|d_xv^-|\leq K$ and $|d_xv^+|\leq K$ for all $x\in\cU$. Since $H$ is Tonelli, the positive definiteness of $D_{pp}^2H$ yields a constant $\lambda>0$ such that
\begin{equation}\label{uniconvex}
	D_{pp}^2H(x,p)\succeq\lambda I_n,\quad \text{for all}~x\in M, \quad |p|\leq K. 
\end{equation}

Fix $x\in\cU$. Set $p^-:=d_xv^-$ and $p^+:=d_xv^+$. Consider
\[
\varphi(t)
:=
H\bigl(x,(1-t)p^-+tp^+\bigr),
\qquad t\in[0,1].
\]
By \eqref{uniconvex}, $\varphi''(t)\geq\lambda\,|p^+-p^-|^2$ for all $t\in [0,1]$. Then the strong convexity inequality yields
\[
\varphi(t)
\leq
(1-t)\varphi(0)+t\varphi(1)
-\frac{\lambda}{2}t(1-t)|p^+-p^-|^2.
\]
Taking $t=1/2$ and using \eqref{Hxvpm_cH} gives
\begin{equation*}
\begin{aligned}
H\left(x,\frac{p^-+p^+}{2}\right)
\leq 
\frac{H(x,p^-)+H(x,p^+)}{2}
-\frac{\lambda}{8}|p^--p^+|^2= c_H-\frac{\lambda}{8}|p^--p^+|^2.
\end{aligned}
\end{equation*}
Combining this with \eqref{eq_dxvpmbeta} and \eqref{defu1}, we obtain
\[
H(x,d_xu_1)\leq~
c_H-\frac{\lambda}{8}|d_xv^--d_xv^+|^2\leq~
c_H-\frac{\lambda\beta^2}{8}d(x,\cA)^2.
\]
Therefore, \eqref{eq_quadraticgap} holds with $a:=\frac{\lambda\beta^2}{8}>0$. This proves the claim.

\medskip

\noindent {\bf Step 2: global extension.} We next extend this local estimate to the whole manifold. By Theorem \ref{thm_smoothCksol}, there exists a $C^k$ critical subsolution $u_0:M\to\R$ which is strict outside $\cA$:
\begin{equation}\label{Bernard_Hu_0}
H(x,d_xu_0)\leq c_H
\quad\text{for all}~x\in M,
\qquad
H(x,d_xu_0)<c_H
\quad\text{for all}~x\in M\setminus\cA.
\end{equation}
Since $u_1$ defined in \eqref{defu1} is $C^k$ on $\cU$, not necessarily on $M$, we introduce a smooth cutoff function supported in $\cU$. More precisely, choose an open neighborhood $O$ of $\cA$ such that $\cA\subset O\subset\bar{O}\subset\cU$, and let $\rho\in C^\infty(M)$ satisfy
\[
0\leq\rho\leq1,
\qquad
\rho=1\ \text{on}~O,
\qquad
\supp\rho\subset \cU.
\]
For $\vep >0$, set
\[
f_\vep(x):=u_0(x)+\vep\,\rho(x)\,[u_1(x)-u_0(x)],
\qquad x\in M.
\]
It follows that $f_\vep: M\to\R$ is of class $C^k$.

On the compact set $M\setminus O$, by \eqref{Bernard_Hu_0}, there exists $\delta>0$ such that $H(x,d_xu_0)\leq c_H-\delta$ on $M\setminus O$. Moreover, $\|f_\vep- u_0\|_{C^1(M)}\to 0$ as $\vep\to 0$, so we can find a small $\vep_0\in (0,1)$ such that 
\[
H(x,d_xf_{\vep_0})\leq c_H-\frac{\delta}{2},
\qquad x\in M\setminus O.
\]
We thus obtain
\begin{equation}\label{quadraticgap_away}
	H(x,d_xf_{\vep_0})\leq c_H-\frac{\delta}{2D^2}\,d(x,\cA)^2,
\qquad x\in M\setminus O,
\end{equation}
where $D:=1+\max_{x\in M}d(x,\cA)$.

 On the set $O$, since $\rho=1$ on $O$, we have
\[
f_{\vep_0}=(1-\vep_0)u_0+\vep_0 u_1 \quad \text{on}~ O.
\]
Using the $p$-convexity of $H$ together with \eqref{eq_quadraticgap} and \eqref{Bernard_Hu_0} yields
\begin{equation}\label{quadraticgap_near}
H(x,d_xf_{\vep_0})
\leq  (1-\vep_0)H(x,d_xu_0)+\vep_0 H(x,d_x u_1)
\leq  c_H-\vep_0 a\,d(x,\cA)^2,
\qquad  x\in O.
\end{equation}

Therefore, combining the above two estimates \eqref{quadraticgap_away} and \eqref{quadraticgap_near}, we obtain
\[
H(x,d_xf_{\vep_0})
\leq
c_H-\kappa\,d(x,\cA)^2,
\qquad x\in M,
\]
where $\kappa:=\min\left\{\vep_0a,\frac{\delta}{2D^2}\right\}>0$. Moreover, since $\rho=1$ near $\cA$, we have
\[ d_xf_{\vep_0}=(1-\vep_0)d_xu_0+\vep_0 d_xu_1,\quad x\in\cA.\]
It is known (see, e.g., \cite{Fathi_Siconolfi2004}) that critical subsolutions have the same differential on the projected Aubry set, so $d_xu_0=d_xu_1$ for all $x\in\cA$. This implies $H(x,d_xf_{\vep_0})=c_H$ for all $x\in\cA$.

Therefore, taking $f:=f_{\vep_0}$ gives the desired $C^k$ critical subsolution. Since the subsolution $f$ is strict outside $\cA$, according to \cite{Fathi_Siconolfi2004}, one has $\tcI^*(f)=\cAt^*$.
\end{proof}

\section{Action Estimates for Connecting Curves Near the Aubry Set}
\label{sec_connections}

We apply Theorem \ref{Thm_quadratic} to derive quantitative action estimates for connecting curves near the Aubry set. These estimates will be used to establish the convergence rates in Section \ref{sec_convergencerate}.

\subsection{The normalized Lagrangian}
To simplify the analysis in the following sections, we introduce a standard normalization of the Lagrangian, similar to that in \cite{Morgado15_hyperbolicity}.

 Fix the subsolution $f\in C^k(M)$ from Theorem \ref{Thm_quadratic}.  We introduce the {\it normalized Lagrangian}
\begin{equation}\label{eq_newL}
\hL(x,v):=L(x,v)-d_xf(v)+c_H,
\qquad (x,v)\in TM.
\end{equation}
Such a modification is standard in the study of Lagrangian dynamics (see, e.g., \cite{Mather1991, Mather1993, Sorrentino_book, ChengXuebook}). Since $df$ is exact, $L$ and $\hL$ generate the same Euler–Lagrange flow on $TM$, even though $\hL$ is generally only of class $C^{k-1}$. Moreover, $\hL$ enjoys the following basic properties.

\begin{lemma}\label{lem_newL_properties}
Let $\hL: TM\to \R$ be defined by \eqref{eq_newL}. Then the following statements hold.
\begin{enumerate}
\item[\rm(i)]
For every continuous piecewise $C^1$ curve $\gamma:[a,b]\to M$, define its $\hL$-action by
\begin{equation*}
\nAct(\gamma):=\int_a^b\hL(\gamma(s),\dot\gamma(s))\,ds.
\end{equation*}
Then
\begin{equation*}
\nAct(\gamma)=A_L(\gamma)+c_H(b-a)+f(\gamma(a))-f(\gamma(b)).
\end{equation*}
Consequently, the minimizing curves of $\hL$ coincide with those of $L$.

\item[\rm(ii)] We have $\hL\geq 0$ on $TM$, and 
\begin{equation}\label{levelset}
\hL(x,v)=0
\quad\Longleftrightarrow\quad
(x,v)\in\cAt,
\end{equation}
where $\cAt$ is the Aubry set of the original Lagrangian $L$. Moreover, 
\begin{equation}\label{eq_newL_lowerbd}
\hL(x,v)\geq \kappa\,d(x,\cA)^2,\quad \forall~(x,v)\in TM,
\end{equation}
where $\kappa>0$ is the constant given in Theorem \ref{Thm_quadratic}.

\item[\rm(iii)] For every nonempty compact subset $\mathcal{K}\subset TM$, there exists $c_{\mathcal{K}}>0$ such that
\begin{equation}\label{eq_quaddist_lbd}
\hL(x,v)\ge c_{\mathcal{K}}\left[d_{TM}\big((x,v),\cAt\big)\right]^2,
\quad \forall~(x,v)\in \mathcal K,
\end{equation}
where $d_{TM}$ is a fixed Riemannian distance on $TM$.
\end{enumerate}
\end{lemma}

\begin{proof}
(i) Since $\hL(x,v)=L(x,v)-d_xf(v)+c_H$, the chain rule gives
\[
\begin{aligned}
\nAct(\gamma)
&=
\int_a^b
\bigl(L(\gamma(s),\dot\gamma(s))+c_H\bigr)\,ds
-\int_a^b d_{\gamma(s)}f(\dot\gamma(s))\,ds
\\
&=
A_L(\gamma)+c_H(b-a)
+f(\gamma(a))-f(\gamma(b)).
\end{aligned}
\]

(ii)  We begin with \eqref{eq_newL_lowerbd}. Indeed, by the Fenchel inequality and Theorem \ref{Thm_quadratic},
\[
\hL(x,v)=L(x,v)-d_xf(v)+c_H\geq c_H-H(x,d_xf)\geq \kappa d(x,\cA)^2.
\]
In particular, $\hL\geq0$ on $TM$. It remains to prove \eqref{levelset}. Write
\[
\hL(x,v)=\big[L(x,v)+H(x,d_xf)-d_xf(v)\big]
+
\big[c_H-H(x,d_xf)\big].
\]
Thus, $\hL(x,v)=0$ if and only if $L(x,v)+H(x,d_xf)=d_xf(v)$ and $H(x,d_xf)=c_H$. Consequently, by the equality condition in Fenchel's inequality and Theorem \ref{Thm_quadratic}, we obtain
\begin{equation}\label{lzeroset}
\hL(x,v)=0
\quad\Longleftrightarrow\quad
v=\partial_pH(x,d_xf)
\quad \text{and}\quad
x\in\cA.
\end{equation}

By Theorem \ref{Thm_quadratic}, $\cAt^*=\tcI^*(f)$. Hence Mather's graph property
(Proposition \ref{prop_Matherproperty}), together with $\pi(\cAt^*)=\cA$, yields
\[
\cAt^*=\{(x,d_xf):x\in\cA\}.
\]
Recalling that $\cAt={\mathscr L}^{-1}(\cAt^*)$ via the Legendre transform $\mathscr L$,  we conclude from \eqref{lzeroset} that $\hL(x,v)=0$ if and only if $(x,v)\in\cAt$.

(iii)
Set $\hat{v}(x):=H_p(x,d_xf)$ on $M$. Since $\partial_{vv}\hL=\partial_{vv}L$ is positive definite, compactness and Taylor's theorem imply that there exists $c_1>0$ such that
\begin{equation*}
	\begin{aligned}
		\hL(x,v)\ge \hL(x,\hat{v}(x))+\partial_v \hL(x,\hat{v}(x))\,(v-\hat{v}(x))+c_1 |v-\hat{v}(x)|_x^2
	\end{aligned}
\end{equation*}
for all $(x,v)\in\mathcal K$.
By Legendre duality, $\partial_v\hL(x,\hat{v}(x))=\partial_v L(x,\hat{v}(x))-d_xf=0$.  Together with \eqref{eq_newL_lowerbd}, this yields
\begin{equation*}
\hL(x,v)\geq \kappa\, d(x,\cA)^2+c_1|v-\hat{v}(x)|_x^2\geq c_2 \big(\,d(x,\cA)+|v-\hat{v}(x)|_x\,\big)^2,\qquad \forall~(x,v)\in\mathcal K
\end{equation*}
for some $c_2>0$. Next, for each fixed $(x,v)\in \mathcal{K}$, we choose $a\in \cA$ such that $d(x, a)=d(x, \cA)$. Since $\cAt=\{(y,\hat{v}(y)):y\in\cA\}$, $\hat{v}$ is at least $C^1$ on $M$, and $\mathcal K$ is compact, we obtain
\[
\begin{aligned}
d_{TM}\big((x,v), \cAt\big)
&\le d_{TM}\big((x,v),(x,\hat{v}(x))\big)
+d_{TM}\big((x, \hat{v}(x)),(a, \hat{v}(a))\big)\\
&\le c_3\big(\,|v-\hat{v}(x)|_x+d(x,a)\,\big)= c_3\big(\,|v-\hat{v}(x)|_x+d(x,\cA)\,\big)
\end{aligned}
\]
for some $c_3>0$. Therefore, \eqref{eq_quaddist_lbd} follows with $c_{\mathcal K}=c_2/c_3^2$. 
\end{proof}

\begin{remark}
 Since $\hL\geq0$ on $TM$ and $\hL=0$ on $\cAt$ by \eqref{levelset}, it follows immediately that
\begin{equation}\label{eq_newL_dv}
\partial_v\hL(x,v)=0,
\quad\text{for all}~ (x,v)\in\cAt.
\end{equation}
The lower bound \eqref{eq_newL_lowerbd} will be used extensively in the subsequent sections.
\end{remark}

\subsection{Connecting curves and loops with controlled action}
In the Lagrangian setting, the Aubry set $\cAt\subset TM$ of $L$ is $\cAt=\cAt_1\cup\cdots\cup\cAt_m,$ where $\cAt_i:=\mathscr L^{-1}(\cAt_i^*)$. Each $\cAt_i$ is either a hyperbolic periodic orbit or a hyperbolic fixed point of the Euler--Lagrange flow. We fix some notation for the constructions below:
\begin{itemize}

\item  For each $i$, let $\bgamma_i:\R\to M$ be a curve whose lift $t\mapsto(\bgamma_i(t),\dbgamma_i(t))$ is a trajectory of the Euler--Lagrange flow with image $\cAt_i$, that is,
\begin{equation}\label{trajectoryimage}
\cAt_i
=
\big\{(\bgamma_i(t),\dbgamma_i(t)):t\in\R\big\}.
\end{equation}
Thus, $\cA_i=\{\bgamma_i(t):t\in\R\}$ is either an embedded $C^k$ circle or a singleton.

\item Let $\mathbf{I}_{\mathrm{per}}$ denote the set of indices corresponding to periodic orbits, and let $\mathbf{I}_{\mathrm{fix}}:=\{1,\ldots,m\}\setminus \mathbf{I}_{\mathrm{per}}$ be the set of indices corresponding to fixed points. For each $i$, define
\begin{equation}\label{convention_bTi}
\bT_i:=
\begin{cases}
\text{the least period of}~\bgamma_i, & \quad i\in\mathbf{I}_{\mathrm{per}}\,,\\
1, & \quad i\in \mathbf{I}_{\mathrm{fix}}\,.
\end{cases}
\end{equation}
We remark that in the fixed-point case, $\cAt_i=\{(a_i,0)\}$, $\bgamma_i(s)\equiv a_i$, and the choice $\bT_i=1$ is merely a convention (not a least period).

\item For $r>0$, we denote the open $r$-neighborhood of $\cA_i$ in $M$ by
\[
V_r(\cA_i)
:=
\{x\in M:d(x,\cA_i)<r\}.
\]
\end{itemize}

\begin{lemma}[Connecting curves near the Aubry set]\label{lem:bridges}
There exist constants $r_0>0$ and $C_{\mathrm b}>0$ such that, for each $i=1,\ldots,m$, the following hold:
\begin{enumerate}
\item[\rm(i)] For any two points $x_s, x_e\in V_{r_0}(\cA_i)$, there exist $\tau\in [2,2+\bT_i]$ and a continuous piecewise $C^1$ curve $\eta:[0,\tau]\to M$ joining $x_s$ to $x_e$, such that $\eta([0,\tau])\cap\cA_i\neq\emptyset$ and
\begin{equation*}
\int_0^\tau \hL(\eta(t),\dot\eta(t)) \,dt
\leq
C_{\mathrm b}
\left[
d(x_s,\cA_i)^2+d(x_e,\cA_i)^2
\right].
\end{equation*}

\item[\rm(ii)] For any $x\in V_{r_0}(\cA_i)$, there exists $\tau\in[2,2+\bT_i]$ such that, for every integer $N\geq 1$, there is a continuous piecewise $C^1$ loop $\eta_N:[0,\tau+N\bT_i]\to M$ based at $x$ satisfying 
\[
\int_0^{\tau+N\bT_i} \hL(\eta_N(t),\dot\eta_N(t)) \,dt
\leq
2C_{\mathrm b}\,d(x,\cA_i)^2.
\]
\end{enumerate}
\end{lemma}

\begin{proof}
{\bf Step 1.} We first construct short connecting curves near each point of $\cA=\cup_{i=1}^m \cA_i$. The construction below applies to both the periodic-orbit case and the fixed-point case.

Fix $q\in\cA_i$. By \eqref{trajectoryimage}, we choose $t_0$ such that $q=\bgamma_i(t_0)$. Let $U$ be a smooth local coordinate chart around $q$. Choose  suitably small $\delta\in(0,1)$ and $r_0>0$ such that
\[
\bgamma_i\big|_{[t_0-\delta,t_0+\delta]}\subset U,
\qquad
B_{r_0}(q)\subset U.
\]
To simplify notation, we identify $U$ with an open subset of $\R^n$ via the coordinate map. 

For every $y\in B_{r_0}(q)$, we define a curve $\xi_y^s:[0,\delta]\to U\subset M$, joining $y$ to $\bgamma_i(t_0+\delta)$, by
\[
\xi_y^s(t)
:=
\bgamma_i(t_0+t)
+\frac{\delta-t}{\delta}(y-q),\qquad t\in[0,\delta].
\]
See Figure \ref{fig_step1shortcurve} for an illustration of the construction.
\begingroup
\setlength{\intextsep}{5pt}
\setlength{\abovecaptionskip}{4pt}
\setlength{\belowcaptionskip}{0pt}
\begin{figure}[H]
  \centering
  \includegraphics[scale=1.2]{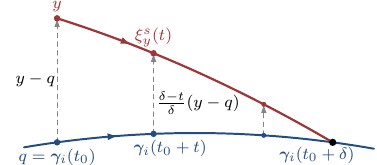}
 \caption{Construction of the connecting curve $\xi_y^s$}
 \label{fig_step1shortcurve}
\end{figure}
\endgroup
\noindent Let $F_s(y):=A_{\hL}(\xi_y^s)$. Since
$
\hL(x,v)=L(x,v)-d_xf(v)+c_H
$, 
\[
\begin{aligned}
F_s(y)
&=\int_0^\delta
L\left(
\bgamma_i(t_0+t)+\frac{\delta-t}{\delta}(y-q),
\dbgamma_i(t_0+t)-\frac{(y-q)}{\delta}
\right)
dt
\\
&\qquad +f(y)-f(\bgamma_i(t_0+\delta))+c_H\delta.
\end{aligned}
\]
It follows from this expression that $F_s$ is $C^k$ $(k\geq 2)$ in $y$, since $L, f\in C^k$. By Lemma \ref{lem_newL_properties}(ii),
\begin{equation}\label{bgamma0action}
\hL\ge 0~\text{on}~TM,
\qquad
\hL\bigl(\bgamma_i(t),\dbgamma_i(t)\bigr)=0
~\text{for all}~t\in\R.
\end{equation}
Consequently, $F_s(y)\geq 0$, $F_s(q)=0$ and $DF_s(q)=0$. Thus, Taylor's theorem gives 
\begin{equation*}
F_s(y)\le C|y-q|^2,
\qquad
\forall~y\in B_{r_0}(q),
\end{equation*}
for some constant $C>0$, after decreasing $r_0$ if necessary.

Similarly,  for every $y\in B_{r_0}(q)$, define a curve $\xi_y^e$ in $U$, connecting $\bgamma_i(t_0-\delta)$ to $y$, by
\[
\xi_y^e(t)
:=
\bgamma_i(t_0+t)+\frac{t+\delta}{\delta}(y-q),
\qquad -\delta\leq t\leq 0.
\]
Let $F_e(y):=A_{\hL}(\xi_y^e)$. The same argument shows that
\begin{equation*}
F_e(y)\leq C|y-q|^2,
\qquad
\forall~y\in B_{r_0}(q).
\end{equation*}

We now make these estimates uniform in $q$. By compactness of $\cA$, choose finitely many coordinate neighborhoods $U_j$ and open subsets $W_j\Subset U_j$ such that $\cA\subset\bigcup_j W_j$. We may choose $\delta\in(0,1)$ and $r_0>0$ uniformly so that the above constructions are valid in $U_j$ whenever $q\in\cA\cap W_j$ and $d(y,q)<r_0$. The compactness yields uniform bounds on the second derivatives of $F_s$ and $F_e$, as well as uniform equivalence of Euclidean and Riemannian distances. Thus, Taylor's theorem gives a constant $C_{\mathrm b}>0$, independent of $q$, such that
\begin{equation}\label{eq_shortbridgeaction}
F_s(y)\leq C_{\mathrm b}\,d(y,q)^2,
\qquad
F_e(y)\leq C_{\mathrm b}\,d(y,q)^2,
\end{equation}
whenever $q\in\cA$ and $d(y,q)<r_0$.

\medskip

\noindent {\bf Step 2.} We now prove assertion (i).\,  Fix $x_s,x_e\in V_{r_0}(\cA_i)$ and choose $q_s,q_e\in\cA_i$ so that
\begin{equation}\label{distxpm}
d(x_s,q_s)=d(x_s,\cA_i),
\qquad
d(x_e,q_e)=d(x_e,\cA_i).
\end{equation}
We may choose $\theta_s, \theta_e\in \R$ such that 
\[
q_s=\bgamma_i(\theta_s),
\qquad
q_e=\bgamma_i(\theta_e),\qquad 0\le \theta_e-\theta_s<\bT_i.
\]
We remark that we may take $\theta_s=\theta_e=0$ when $i\in \mathbf{I}_{\mathrm{fix}}$ (i.e.,  the fixed-point case).

Applying the construction in Step 1 separately near $q_s$ and $q_e$, we obtain two curves
\[
\xi_s:[0,\delta]\to M,
\qquad
\xi_e:[-\delta,0]\to M
\]
satisfying
\begin{equation}\label{endpts}
\xi_s(0)=x_s,
\quad
\xi_s(\delta)=\bgamma_i(\theta_s+\delta);
\qquad\quad
\xi_e(-\delta)=\bgamma_i(\theta_e-\delta),
\quad
\xi_e(0)=x_e.
\end{equation}
Applying \eqref{eq_shortbridgeaction} and \eqref{distxpm}, we obtain
\begin{equation}\label{twocurves_action}
\int_0^\delta
\hL(\xi_s(t),\dot \xi_s(t))\,dt
+\int_{-\delta}^0
\hL(\xi_e(t),\dot \xi_e(t))\,dt
\leq
C_{\mathrm b}\,\big[d(x_s,\cA_i)^2+d(x_e,\cA_i)^2\big].
\end{equation}

Choose $\tau\in[2,2+\bT_i]$ such that
\begin{equation}\label{tauchoice}
\tau-(\theta_e-\theta_s)\in\bT_i\mathbb Z.
\end{equation}
Such a choice is possible because the interval $[2,2+\bT_i]$ has length $\bT_i$. We also remark that when $i\in \mathbf{I}_{\mathrm{fix}}$, our convention (see \eqref{convention_bTi}) $\bT_i=1$ allows us to choose $\tau=2$. Since $\delta<1$, one has $2\delta<\tau$. Define $\eta: [0,\tau]\to M$ as the concatenation of three curve segments:  
\begin{equation}\label{curve_concatenation}
\eta(t):=
\begin{cases}
\xi_s(t),
& \quad t\in[0,\delta],
\\
\bgamma_i(\theta_s+t),
& \quad t\in[\delta,\tau-\delta],
\\
\xi_e(t-\tau),
& \quad t\in [\tau-\delta,\tau].
\end{cases}
\end{equation}
By \eqref{endpts} and \eqref{tauchoice}, $\eta$ is a continuous piecewise $C^1$ curve satisfying
\[
\eta(0)=x_s,
\qquad
\eta(\tau)=x_e.
\]
Moreover, its middle segment $\eta\big|_{[\delta, \tau-\delta]}$ has zero $\hL$-action by \eqref{bgamma0action}. This together with \eqref{twocurves_action} proves assertion (i).

\medskip
\noindent{\bf Step 3.}
It remains to prove assertion (ii). Applying assertion (i) with $x_s=x_e=x\in V_{r_0}(\cA_i)$, we obtain a continuous piecewise $C^1$ loop $\eta:[0,\tau]\to M$ satisfying $\eta(0)=\eta(\tau)=x$ and
\[
\int_0^\tau
\hL(\eta(t),\dot\eta(t))\,dt
\leq 2C_{\mathrm b}\,d(x,\cA_i)^2.
\]
By \eqref{curve_concatenation}, there exist $\sigma\in(0,\tau)$ and $t_0\in [0,\bT_i)$ such that $\eta(\sigma)=\bgamma_i(t_0)\in \cA_i$. 

For each $N\geq1$, the trajectory segment $\bgamma_i|_{[t_0,t_0+N\bT_i]}$ is a loop based at $\eta(\sigma)$. Inserting this loop into $\eta$ at time $\sigma$ yields a continuous piecewise $C^1$ loop $\eta_N:[0,\tau+N\bT_i]\to M$ satisfying $\eta_N(0)=\eta_N(\tau+N\bT_i)=x$. Moreover, by \eqref{bgamma0action}, we deduce that
\[
\begin{aligned}
\int_0^{\tau+N\bT_i}
\hL(\eta_N(t),\dot\eta_N(t))\,dt
&=
\int_0^\tau
\hL(\eta(t),\dot\eta(t))\,dt
\leq
2C_{\mathrm b}\,d(x,\cA_i)^2.
\end{aligned}
\]
This proves assertion (ii).
\end{proof}

\begin{lemma}[Loops of prescribed duration]
\label{lem_loops}
Set $\bT_*:=\max_{1\leq i\leq m}\bT_i$. There exists a constant $C_{\mathrm p}>0$ such that,
for every $x\in\cA$ and every $t\geq2\bT_*$,
there exists a $C^1$ loop $\zeta:[0,t]\to M$
satisfying $\zeta(0)=\zeta(t)=x$ and 
\begin{equation}\label{eq_loopsaction}
\int_0^t \hL(\zeta(s),\dot\zeta(s)) \,ds \leq \frac{C_{\mathrm p}}{t}.
\end{equation}
\end{lemma}

\begin{proof}
Fix $x\in\cA$ and $t\geq2\bT_*$, and choose $i\in\{1,\ldots,m\}$ such that $x\in \cA_i$. If $i\in\mathbf{I}_{\mathrm{fix}}$, then $\bgamma_i(s)\equiv a_i$ for all $s\in\R$ and hence $x=a_i$. Taking the constant loop $\zeta(s)\equiv a_i$ on $[0,t]$, we have $\int_0^t\hL(\zeta(s),\dot\zeta(s))\,ds=t\hL(a_i,0)=0$. Thus \eqref{eq_loopsaction} holds trivially.

Now suppose that $i\in\mathbf{I}_{\mathrm{per}}$. Choose $s_0\in\R$ such that $x=\bgamma_i(s_0)$. Since $t\geq 2\bT_*$, we may choose an integer $N\geq 1$ such that $\left|t-N\bT_i\right|\leq\frac{\bT_i}{2}$. Set $\lambda:=N\bT_i/t$. Then
\begin{equation}\label{NTiest}
|\lambda-1|\leq\frac{\bT_i}{2t}\leq\frac14,
\qquad \text{and hence}\qquad
\lambda\in\left[\frac34,\frac54\right].
\end{equation}
Define 
\[
\zeta(s)
:=
\bgamma_i\left(
s_0+\lambda s
\right),
\qquad s\in [0,t].
\]
By the $\bT_i$-periodicity of $\bgamma_i$, $\zeta$ is a $C^1$ loop satisfying $\zeta(0)=\zeta(t)=x$.

By Lemma \ref{lem_newL_properties}(ii) and \eqref{eq_newL_dv}, we have $\hL\bigl(\bgamma_j(\tau),\dbgamma_j(\tau)\bigr)=0$ and $\partial_v\hL\bigl(\bgamma_j(\tau),\dbgamma_j(\tau)\bigr)=0$ for all $\tau\in\R$. Then Taylor's theorem yields a constant $C_0>0$ such that
\[
\hL\bigl(\bgamma_j(\tau),\lambda\,\dbgamma_j(\tau)\bigr)
\leq
C_0\,|\lambda-1|^2,\qquad \text{for all $\tau\in \R$, $j\in\mathbf{I}_{\mathrm{per}}$, $\lambda\in\left[\frac34,\frac54\right]$}.
\]
Combining this with \eqref{NTiest} gives
\[
\begin{aligned}
\nAct(\zeta)
=
\int_0^t
\hL\left(
\bgamma_i\left(s_0+\lambda s\right),
\lambda\dbgamma_i\left(s_0+\lambda s\right)
\right)\,ds \leq
C_0\,t\,
\left|
\lambda-1
\right|^2 \leq\frac{C_0\bT_i^2}{4t}\,\leq\,\frac{C_0\bT_*^2}{4t}.
\end{aligned}
\]
Thus \eqref{eq_loopsaction} holds with $C_{\mathrm p}:=\frac{C_0\bT_*^2}{4}$.
\end{proof}

\section{Convergence Rate Estimates}
\label{sec_convergencerate}

Let $\hL: TM\to\R$ be the modified Lagrangian given in \eqref{eq_newL}. For $t>0$, let $\kt:M\times M\to \R$ denote the time-$t$ minimal action of $\hL$, defined by
\begin{equation*}
\kt(y,x):=\min_{\xi\in\Gamma^t_{y,x}}\nAct(\xi).
\end{equation*}
Let $\kinfty:M\times M\to \R$ denote the corresponding Peierls barrier:
\[
\kinfty(y,x)
:=\liminf_{t\to\infty}\kt(y,x).
\]
Since $\hL(x,v)=L(x,v)-d_xf(v)+c_H$, it follows easily from Lemma \ref{lem_newL_properties}(i) that 
\begin{equation}\label{relation_ktht}
\begin{aligned}
\kt(y,x)&=h_t(y,x)+c_Ht-f(x)+f(y),\\
\kinfty(y,x)&=h(y,x)-f(x)+f(y),
\end{aligned}
\end{equation}
where $h_t$ and $h$ are, respectively, the minimal action and the Peierls barrier of the original Lagrangian $L$. It follows immediately that
\begin{equation}\label{diff_ktkinfty}
	\kt(y,x)-\kinfty(y,x)=h_t(y,x)+c_Ht-h(y,x).
\end{equation}

\begin{theorem}\label{thm_equivstatement}
There exist constants $C>0$ and $t_*>0$ such that, for every $t\geq t_*$,
\begin{equation*}
	\sup_{x,y\in M}| h_t(y,x)+c_Ht-h(y,x)|=\sup_{x,y\in M}| \kt(y,x)-\kinfty(y,x)|\leq C t^{-1}.
\end{equation*}
\end{theorem}

Theorem \ref{thm_equivstatement} thus yields the upper-bound estimate \eqref{h_upperbd} in Theorem \ref{Thm_kernel}. The lower-bound estimate in \eqref{h_lowerbd} will be established at the end of Section \ref{sec_optimal} as a consequence of Theorem \ref{Thm_optimality}.

\subsection{Proof of Theorem \ref{thm_equivstatement}}
We begin with the following basic estimates.

\begin{lemma}\label{lem_Erestimate}
There exists a constant $B>0$ such that
\begin{equation}\label{ktB}
0\leq \kt(y,x)\leq B \qquad\text{for all }x,y\in M,\ t\ge1.
\end{equation}
Moreover, for every $x,y\in M$ and $t\ge1$, any minimizer $\gamma:[0,t]\to M$ for $\kt(y,x)$ satisfies 
\begin{equation}\label{eq_distanceintegral}
\int_0^t d(\gamma(s),\cA)^2\,ds
\le B\kappa^{-1},
\end{equation}
where $\kappa>0$ is the constant given in Theorem \ref{Thm_quadratic}. Consequently, for any $r>0$, setting
$
E_r:=\{s\in[0,t]:d(\gamma(s),\cA)>r\},
$
we have
\begin{equation}\label{eq_occupation}
\meas\left(E_r\right)
\le (B \kappa^{-1})\, r^{-2},
\end{equation}
where $\meas$ denotes the Lebesgue measure. 
\end{lemma}

\begin{remark}
Estimate \eqref{eq_occupation} shows that the total time $\gamma$ spends outside the closed $r$-neighborhood $\overline{V_r(\cA)}$ is at most $B/(\kappa r^2)$.
\end{remark}

\begin{proof}
Since $\hL\geq 0$, we clearly have $\kt(y,x)\geq 0$. Meanwhile, \eqref{relation_ktht} gives 
\[
\|\kt\|_\infty\leq \|h_t+c_Ht\|_\infty+2\|f\|_\infty.
\]
By Proposition \ref{Prop_Fleming},  $\|h_t+c_Ht\|_\infty$ is uniformly bounded for  $t\geq 1$. Hence \eqref{ktB} follows. Moreover, using \eqref{eq_newL_lowerbd} from Lemma \ref{lem_newL_properties},
\[
\kappa\int_0^t d(\gamma(s),\cA)^2\,ds
\leq
\int_0^t\hL(\gamma(s),\dot\gamma(s))\,ds
=
\kt(y,x)
\leq B,
\]
which proves \eqref{eq_distanceintegral}. In particular, \eqref{eq_occupation} follows easily from the following inequality 
\[
r^2\mathcal L^1(E_r)
\le
\int_{E_r}d(\gamma(s),\cA)^2\,ds
\le
B/\kappa.
\]
\end{proof}

Since the sets $\cA_1,\ldots,\cA_m$ are compact and pairwise disjoint, we may decrease the radius $r_0$ in Lemma \ref{lem:bridges}, if necessary, so that their closed $r_0$-neighborhoods $\overline{V_{r_0}(\cA_1)},\ldots, \overline{V_{r_0}(\cA_m)}$ are pairwise disjoint. With $r_0$ thus fixed and $B$ as in Lemma \ref{lem_Erestimate}, we fix $t_*>0$ satisfying
\begin{equation}\label{choice_tstar}
t_*>\frac{4B}{\kappa \, r_0^2},
\qquad
t_*\geq 4m(2+\bT_*),
\end{equation}
where $\bT_*:=\max_{1\leq i\leq m}\bT_i$ with $\bT_i>0$ defined in \eqref{convention_bTi}.

\begin{proposition}[The upper estimate]\label{prop_ktupper}
Let $t_*$ be as in \eqref{choice_tstar}. Then there exists a constant $C_1>0$ such that, for every $x,y\in M$ and $t\geq t_*$,
\begin{equation}\label{eq_ktupper}
\kt(y,x)\le \kT(y,x)
+C_1 t^{-1}, \quad \text{for all}~T\in [t,\infty).
\end{equation}
Consequently, for every $x,y\in M$ and $t\geq t_*$,
\[
\kt(y,x)-\kinfty(y,x)\leq C_1 t^{-1}.
\]
\end{proposition}

\begin{proof}
Fix $x,y\in M$ and $t\geq t_*$. For this $t$, set
\begin{equation}\label{r_choice}
	r:=\sqrt{\frac{4B}{\kappa\, t}}\qquad\text{and}\qquad  \cA_i(r):=\overline{V_r(\cA_i)}=\{z\in M:d(z,\cA_i)\leq r\}.
\end{equation}
Then, $r<r_0$ as a result of $t\geq t_*$ and \eqref{choice_tstar}. Hence, the closed $r$-neighborhoods $\cA_i(r)$, $ i=1,\ldots,m$, are pairwise disjoint. 

Let $T\in [t,\infty)$, and let $\gamma:[0,T]\to M$ be a minimizing curve for $\kT(y,x)$, namely,
\[
\gamma(0)=y,
\qquad
\gamma(T)=x,
\qquad
\nAct(\gamma)=\kT(y,x).
\]
Starting from $\gamma$, we will construct below a curve $\tilde\gamma$ joining $y$ to $x$ in time $t$ and satisfying
\[
\nAct(\tilde\gamma)\leq \kT(y,x)+C_1 t^{-1}.
\]
Our construction is similar to those used in \cite{Iturriaga_Morgado09, Wang_Yan12CMP, Wang_Yan2012nonli, Morgado15_hyperbolicity}. We proceed in three steps.

\medskip

\noindent\textbf{Step 1.} Applying \eqref{eq_occupation} from Lemma \ref{lem_Erestimate}, together with \eqref{r_choice}, we obtain
\begin{equation}\label{meast4}
\meas\left\{
s\in[0,T]:\,
\gamma(s)\notin
\bigcup_{i=1}^m \cA_i(r)
\right\}
\leq
\frac{B}{\kappa \, r^2}
=
\frac{t}{4}.
\end{equation}
This implies that $\gamma:[0, T]\to M$ must meet $\bigcup_{i=1}^m \cA_i(r)$. We now select successively the first and last visits of $\gamma$ to these neighborhoods. 
\begin{enumerate}
	\item [i)] Let
\[
s_1
:=
\min\left\{
s\in[0,T]:
\gamma(s)\in
\bigcup_{i=1}^m \cA_i(r)
\right\},
\]
let $\bm{i_1}\in \{1,\ldots,m\}$ be the corresponding index such that
$
\gamma(s_1)\in \cA_{\bm{i_1}}(r)
$,
and define
\[
e_1
:=
\max\left\{
s\in[s_1,T]:
\gamma(s)\in \cA_{\bm{i_1}}(r)
\right\}.
\]
More intuitively, $s_1$ is the first visit of $\gamma$ to these neighborhoods, $\cA_{\bm{i_1}}(r)$ is the neighborhood containing $\gamma(s_1)$, and $e_1$ is its last visit to $\cA_{\bm{i_1}}(r)$. %In this sense, $s_1$ and $e_1$ means the ``start'' and ``end'' times. 
Notice that $\gamma$ need not remain in $\cA_{\bm{i_1}}(r)$ throughout $[s_1,e_1]$. It may leave and return several times, but it never visits $\cA_{\bm{i_1}}(r)$ again after $e_1$.

\item [ii)] Proceeding inductively, suppose that distinct indices $\bm{i_1},\ldots,\bm{i_j}$ and intervals $[s_1,e_1]$, $\ldots$, $[s_j,e_j]$ have been selected. If $\gamma$ does not meet any of the remaining closed $r$-neighborhoods after $e_j$, the procedure stops. Otherwise, set
\[
s_{j+1}
:=
\min\left\{
s\in[e_j,T]:
\gamma(s)\in
\bigcup_{i\notin\{\bm{i_1},\ldots,\bm{i_j}\}}
\cA_i(r)
\right\}.
\]
Let $\bm{i_{j+1}}$ be the corresponding index such that
$
\gamma(s_{j+1})
\in
\cA_{\bm{i_{j+1}}}(r)
$,
and define
\[
e_{j+1}
:=
\max\left\{
s\in[s_{j+1},T]:
\gamma(s)\in
\cA_{\bm{i_{j+1}}}(r)
\right\}.
\]
\end{enumerate}
This procedure stops after some $\bm{m}_0\leq m$ steps, yielding distinct indices $\bm{i_1},\ldots,\bm{i_{m_0}}$ and disjoint intervals
\[
[s_1,e_1], [s_2,e_2],\ldots, [s_{\bm{m}_0},e_{\bm{m}_0}].
\]
Moreover, the total length of the remaining time intervals satisfies
\begin{equation}\label{remaining_duration}
\sigma:=T-\sum_{j=1}^{\bm{m}_0}(e_j-s_j)\leq\frac{t}{4}.
\end{equation}
Indeed, $
[0,T]\setminus\bigcup_{j=1}^{\bm{m}_0}[s_j,e_j]
\subset
\left\{
s\in[0,T]:
\gamma(s)\notin
\bigcup_{i=1}^m \cA_i(r)
\right\},
$
so \eqref{remaining_duration} follows from \eqref{meast4}.

\medskip

\noindent\textbf{Step 2.}
Next, we replace each curve segment $\gamma|_{[s_j,e_j]}$ by a curve of uniformly bounded duration and small action.

For each $j=1,\ldots,\bm{m}_0$, since
$\gamma(s_j),\gamma(e_j)
\in
\cA_{\bm{i_j}}(r)
\subset
V_{r_0}(\cA_{\bm{i_j}})$,
applying Lemma \ref{lem:bridges}(i) gives a continuous piecewise $C^1$ curve
\[
\eta_j:[0,\tau_j]\to M\quad \text{with}~\tau_j\leq 2+\bT_*,
\]
joining $\gamma(s_j)$ to $\gamma(e_j)$ and satisfying
\begin{equation}\label{etacrossA}
\eta_j([0,\tau_j])\cap \cA_{\bm{i_j}}\neq\emptyset, 	
\end{equation}
and
\begin{equation}\label{etajAction}
\nAct(\eta_j)
\leq
C_{\mathrm b}
\left[
d(\gamma(s_j),\cA_{\bm{i_j}})^2
+
d(\gamma(e_j),\cA_{\bm{i_j}})^2
\right]
\leq
2C_{\mathrm b}r^2.
\end{equation}

We now modify the curve $\gamma$ as follows. For each $j=1,\ldots,\bm{m}_0$, we simply remove the portion $\gamma|_{[s_j,e_j]}$ and replace it by $\eta_j$. Since $\eta_j$ has the same endpoints as $\gamma|_{[s_j,e_j]}$, it can be inserted exactly in place of the removed portion $\gamma|_{[s_j,e_j]}$. In this way, we obtain a new continuous piecewise $C^1$ curve joining $y$ to $x$ (see Figure \ref{fig_concatenation}),
\[
\hat\gamma:[0,\ell]\to M,
\qquad 
\hat\gamma(0)=y, ~\hat\gamma(\ell)=x,
\quad\text{and}~
\ell=\sigma+\sum_{j=1}^{\bm{m}_0}\tau_j.
\]

By \eqref{etajAction} and $\hL\geq0$, the $\hL$-action along $\hat\gamma$ satisfies
\begin{equation}\label{hatgamma_action}
\begin{aligned}
\nAct(\hat\gamma)
=\nAct(\gamma)
-\sum_{j=1}^{\bm{m}_0}
\int_{s_j}^{e_j}\hL(\gamma(s),\dot\gamma(s))\,ds+\sum_{j=1}^{\bm{m}_0}\nAct(\eta_j)
&\leq 
\nAct(\gamma)
+
\sum_{j=1}^{\bm{m}_0}\nAct(\eta_j)\\
&\leq 
\kT(y,x)+2\bm{m}_0C_{\mathrm b}r^2\\
&\leq 
\kT(y,x)
+\frac{8mC_{\mathrm b}B}{\kappa\, t} ,
\end{aligned}
\end{equation}
where the last line follows from \eqref{r_choice} and $\bm{m}_0\leq m$. Moreover, by \eqref{remaining_duration}, \eqref{choice_tstar} and $t\geq t_*$, the total duration $\ell$ of $\hat\gamma$ satisfies
\begin{equation}\label{length_ell}
\ell
\leq
\frac{t}{4}
+
m(2+\bT_*)
\leq \frac{t}{4}+\frac{t_*}{4}\leq 
\frac{t}{2}.
\end{equation}

\begingroup
\setlength{\intextsep}{4pt}
\setlength{\abovecaptionskip}{3pt}
\setlength{\belowcaptionskip}{0pt}
\begin{figure}[H]
  \centering
  \includegraphics[scale=1]{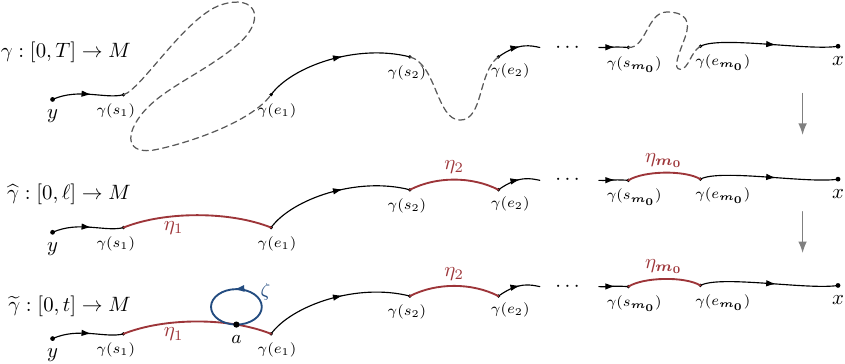}
  \caption{Construction of $\widehat\gamma$ and $\widetilde\gamma$.}
  \label{fig_concatenation}
\end{figure}
\endgroup

\noindent\textbf{Step 3.} We now use Lemma \ref{lem_loops} to extend $\hat\gamma:[0,\ell]\to M$ to a curve of duration $t$ by inserting a loop through the Aubry set.

By \eqref{etacrossA}, each replacement curve $\eta_j$ passes through $\cA_{\bm{i_j}}\subset\cA$. Hence, the concatenated curve $\hat\gamma:[0,\ell]\to M$ also passes through $\cA$, so there exist $\theta\in[0,\ell]$ and $a\in\cA$ such that $\hat\gamma(\theta)=a$. Moreover, \eqref{length_ell} gives $t-\ell\geq\frac{t}{2}\geq \frac{t_*}{2}$, which, together with \eqref{choice_tstar}, yields $t-\ell\geq2\bT_*$. Therefore, Lemma \ref{lem_loops} yields a $C^1$ loop
\[
\zeta:[0,t-\ell]\to M,
\qquad
\zeta(0)=\zeta(t-\ell)=a,
\]
whose $\hL$-action satisfies 
\begin{equation}\label{zeta_action}
\nAct(\zeta)
\leq
\frac{C_{\mathrm p}}{t-\ell}
\leq
\frac{2C_{\mathrm p}}{t}.
\end{equation}
Then, we insert the loop $\zeta$ into $\hat\gamma$ at $a=\hat\gamma(\theta)$, see Figure \ref{fig_concatenation}. More precisely, define $\widetilde\gamma:[0,t]\to M$ by
\[
\widetilde\gamma(s)
:=
\begin{cases}
\hat\gamma(s),
&\quad  s\in [0,\theta],
\\
\zeta(s-\theta),
&\quad s\in [\theta,\theta+t-\ell],
\\
\hat\gamma(s-(t-\ell)),
& \quad s\in[\theta+t-\ell,t].
\end{cases}
\]
Then $\widetilde\gamma$ is a continuous piecewise $C^1$ curve joining
$y$ to $x$ in time $t$. Using \eqref{hatgamma_action} and \eqref{zeta_action} yields
\[
\nAct(\widetilde\gamma)
=\nAct(\hat\gamma)+\nAct(\zeta)\leq
\kT(y,x)
+\frac{8mC_{\mathrm b}B}{\kappa t}
+\frac{2C_{\mathrm p}}{t}.
\]
Therefore, the desired inequality \eqref{eq_ktupper} follows with $C_1:=\frac{8mC_{\mathrm b}B}{\kappa}+2C_{\mathrm p}$. 

Since $T\geq t$ is arbitrary, taking the lower limit as $T\to+\infty$ yields $\kt(y,x)\leq \kinfty(y,x)+C_1/t$ for all $t\geq t_*$. This completes the proof.
\end{proof}

\begin{proposition}[The lower estimate]\label{prop_ktlower}
Let $t_*>0$ be as in \eqref{choice_tstar}. Then there exists a constant $C_2>0$ such that, for all $x,y\in M$ and $t\geq t_*$,
\begin{equation}\label{eq_ktlower}
\kt(y,x)-\kinfty(y,x)\ge-C_2 t^{-1}.
\end{equation}
\end{proposition}

\begin{proof}
Let
$\xi:[0,t]\to M$ be a minimizing curve for $\kt(y,x)$. By Lemma \ref{lem_Erestimate},
\[
\int_0^t d(\xi(s),\cA)^2\,ds\leq \frac{B}{\kappa}.
\]
Hence there exists $s_0\in[0,t]$ such that
\[
x_0:=\xi(s_0),
\qquad
d(x_0,\cA)^2
\leq
\frac{B}{\kappa\, t}.
\]
Then, for some $i\in \{1,\ldots,m\}$, we have
\begin{equation}\label{x0Ai}
d(x_0, \cA_i)=d(x_0,\cA)\leq  \sqrt{\frac{B}{\kappa \, t}}\,.
\end{equation} 
Using $t\geq t_*$ and \eqref{choice_tstar} yields $d(x_0, \cA_i)\leq \sqrt{B/(\kappa \, t_*)}<r_0$. Thus Lemma \ref{lem:bridges}(ii) applies.

Applying Lemma \ref{lem:bridges}(ii), we obtain $\tau\in[2,2+\bT_*]$ such that, for every integer $N\geq1$, there is a continuous piecewise $C^1$ loop $\eta_N:[0,\ell_N]\to M$, where $\ell_N:=\tau+N\bT_i$, satisfying
\[
\eta_N(0)=\eta_N(\ell_N)=x_0
,\qquad 
A_\hL\left(\eta_N\big|_{[0, \ell_N]}\right)
\leq
2C_{\mathrm b}\,d(x_0,\cA_i)^2.
\]

We now concatenate the curve $\xi$ with the loop $\eta_N$ at the point $x_0=\xi(s_0)$. More precisely, define $\bm{\xi_N}:[0,t+\ell_N]\to M$ by
\[
\bm{\xi_N}(s):=
\begin{cases}
\xi(s),
&\quad s\in[0,s_0],
\\
\eta_N(s-s_0),
&\quad s\in [s_0, s_0+\ell_N],
\\
\xi(s-\ell_N),
&\quad s\in [s_0+\ell_N, t+\ell_N].
\end{cases}
\]
Clearly, $\bm{\xi_N}$ is a continuous piecewise $C^1$ curve joining $y$ to $x$. Its action satisfies
\begin{align*}
A_\hL\left(\bm{\xi_N}\big|_{[0, t+\ell_N]}\right)
= 
A_\hL\left(\xi\big|_{[0, t]}\right)
+
A_\hL\left(\eta_N\big|_{[0, \ell_N]}\right)\leq 
\kt(y,x)+2C_{\mathrm b}\,d(x_0,\cA_i)^2.
\end{align*}
Therefore,
\[
h_{t+\ell_N}^{\hL}(y,x)\leq A_\hL\left(\bm{\xi_N}\big|_{[0, t+\ell_N]}\right)
\leq
\kt(y,x)+2C_{\mathrm b}\,d(x_0,\cA_i)^2.
\]

Since $\ell_N\to+\infty$ as $N\to\infty$, by the definition of $\kinfty$ and  \eqref{x0Ai}, we obtain
\[
\begin{aligned}
\kinfty(y,x)
\leq
\liminf_{N\to\infty}h_{t+\ell_N}^\hL(y,x)
&\leq
\kt(y,x)+2C_{\mathrm b}\,d(x_0,\cA_i)^2
\\
&\leq\kt(y,x)+\frac{2C_{\mathrm b}B}{\kappa\,t}.
\end{aligned}
\]
The desired inequality \eqref{eq_ktlower} therefore follows with $C_2:=2C_{\mathrm b}B/\kappa$.
\end{proof}

\smallskip

Finally, combining Propositions \ref{prop_ktupper} and \ref{prop_ktlower} gives
\[
\sup_{x,y\in M}|\kt(y,x)-\kinfty(y,x)|
\le C t^{-1},
\qquad t\ge t_*,
\]
where $C:=\max\{C_1, C_2\}$. This, together with \eqref{diff_ktkinfty}, completes the proof of Theorem \ref{thm_equivstatement}.

\subsection{Convergence rate estimate for the Lax--Oleinik semigroup}
We now prove Theorem \ref{Thm_semigroup} by using Theorem \ref{thm_equivstatement}.

\begin{proof}
We show that for every $u\in C(M)$ and $t>0$,
\begin{equation}\label{ineq_Tleqh}
\|T_t^-u+c_Ht-u_\infty\|_{C(M)} \leq \, \|h_t+c_Ht-h\|_{C(M\times M)}.	
\end{equation}

Fix $u\in C(M)$. By the variational representation of the Lax--Oleinik semigroup, 
\begin{equation}\label{eq_Ttuminus}T_t^-u(x)+c_Ht
=\inf_{y\in M}\{u(y)+h_t(y,x)+c_Ht\},
\qquad x\in M.
\end{equation}
Moreover, by Theorem \ref{uinftyexpression},
\begin{equation}\label{eq_Tuinfty}
u_\infty(x)
=\lim_{t\to+\infty}(T_t^-u(x)+c_Ht)=\inf_{y\in M}\{u(y)+h(y,x)\}.
\end{equation}

Set $E(t):=\|h_t+c_Ht-h\|_{C(M\times M)}$. By the definition of $E(t)$, for every $x,y\in M$,
\[
u(y)+h(y,x)-E(t)
\le u(y)+h_t(y,x)+c_Ht
\le u(y)+h(y,x)+E(t).
\]
Taking the infimum over $y\in M$ and using \eqref{eq_Ttuminus}--\eqref{eq_Tuinfty}, we obtain
\[
u_\infty(x)-E(t)
\le T_t^-u(x)+c_Ht
\le u_\infty(x)+E(t),
\qquad x\in M.
\]
By taking the supremum over $x\in M$, we prove the desired inequality \eqref{ineq_Tleqh}. 

Consequently, Theorem \ref{Thm_semigroup} follows immediately from \eqref{ineq_Tleqh} and Theorem \ref{thm_equivstatement}.
\end{proof}

\section{Optimality of the Convergence Rate}\label{sec_optimal}
In this section, we prove Theorem \ref{Thm_optimality}. By assumption, $\cAt$ contains at least one periodic orbit. Fix such an orbit, denoted by $\cAt_0$, and let $T>0$ be its least period. We can thus write
\[
\cAt_0
=\left\{(\bgamma(s),\dot{\bgamma}(s)):s\in\R\right\}, \qquad \bgamma(s+T)=\bgamma(s),
\]
where $s\mapsto(\bgamma(s),\dot{\bgamma}(s))$ is a periodic trajectory of the Euler--Lagrange flow $\Phi_t^L$. Let $\cA_0=\{\bgamma(s):s\in\R\}$ be its projection onto $M$.  

By Mather's graph property, $\bgamma$ is injective on $[0,T)$, and $\cA_0$ is an embedded $C^k$ circle. In particular, $\bgamma(0)\neq \bgamma(T/2)$. Choose $g\in C^k(M)$ satisfying $0\leq g\leq 1$ and 
\begin{flalign}\label{g_choice}
&\begin{aligned}
\quad\bullet\quad &g(\bgamma(0))=0 \quad\text{and}\quad g(\bgamma(T/2))=1\,;\\
\quad\bullet\quad &g=0\quad\text{in a neighborhood of }\cA\setminus\cA_0.
\end{aligned} &&
\end{flalign}
Consequently, $g$ satisfies the following time-$T$ invariance property:
\begin{equation}\label{eq_ginvariance} 
(g\circ\pi)\big(\Phi_T^L(a)\big)=(g\circ\pi)(a),
\qquad \forall~a\in\cAt.
\end{equation}
Recall that $T>0$ is the least period of $\bgamma$ and $\pi:TM\to M$ is the canonical projection.

\begin{lemma}\label{lem_one_period_lowerbound}
Let $\hL(x,v)=L(x,v)-d_xf(v)+c_H$ be the modified Lagrangian defined in \eqref{eq_newL}. There exists $K\ge1$ such that
\begin{equation*}
h^\hL_T(y,x)\ge\frac1K|g(x)-g(y)|^2,
\qquad \text{for all}~x,y\in M.
\end{equation*}
\end{lemma}

\begin{proof}
Fix $x,y\in M$, and let $\eta:[0,T]\to M$ be a minimizing curve for $h_T^{\hL}(y,x)$. By Lemma \ref{lem_newL_properties}(i), $\eta$ is also a minimizing curve for $L$. The standard {\it a priori} compactness of Tonelli minimizers yields a compact set $\mathcal K\subset TM$, independent of $x$ and $y$, such that $(\eta(s),\dot\eta(s))\in\mathcal K$ for all $s\in[0,T]$. Enlarging $\mathcal K$ if necessary, we may assume that $\cAt\subset\mathcal K$.

Since $g\in C^k(M)$ and the flow $\Phi_t^L$ is of class $C^{k-1}$, by the classical smooth dependence of ODE solutions on initial conditions, there exists a constant $C>0$ such that
\begin{equation}\label{eq_fkr3}
\left|(g\circ\pi)(\Phi_t^L(z))-(g\circ\pi)(\Phi_t^L(z'))\right|
\le C\,d_{TM}(z,z'),
\qquad z,z'\in\mathcal K,\quad |t|\le T.
\end{equation}

Write $z(s):=(\eta(s),\dot\eta(s))$ for $s\in[0,T]$. Since $\eta$ is a
minimizing curve for $L$, we have
\[
y=\pi\Phi_{-s}^L(z(s)),
\qquad
x=\pi\Phi_{T-s}^L(z(s)).
\]
Moreover, for every $z_*\in\cAt$, the flow-invariance of $\cAt$ and \eqref{eq_ginvariance} imply that
\[
g(\pi\Phi_{T-s}^L(z_*))=g(\pi\Phi_{-s}^L(z_*)),\qquad s\in [0,T].
\]
Consequently,
\[
\begin{aligned}
|g(x)-g(y)|
&\le
\left|g(\pi\Phi_{T-s}^L(z(s)))
      -g(\pi\Phi_{T-s}^L(z_*))\right|+
\left|g(\pi\Phi_{-s}^L(z_*))
      -g(\pi\Phi_{-s}^L(z(s)))\right|\\
&\le 2C\,d_{TM}(z(s),z_*),
\end{aligned}
\]
where the last line follows from \eqref{eq_fkr3}. Taking the infimum over $z_*\in\cAt$ gives
\[
|g(x)-g(y)|
\le 2C\,d_{TM}(z(s),\cAt),
\quad \text{for all}~s\in[0,T].
\]
Then, using Lemma \ref{lem_newL_properties}(iii), we obtain
\[
\begin{aligned}
|g(x)-g(y)|^2
\le \frac{4C^2}{T}\int_0^T
d_{TM}(z(s),\cAt)^2\,ds
&\le \frac{4C^2}{Tc_{\mathcal K}}
\int_0^T\hL(\eta(s),\dot\eta(s))\,ds\\
&=\frac{4C^2}{Tc_{\mathcal K}}h_T^{\hL}(y,x).
\end{aligned}
\]
This proves the lemma.
\end{proof}

\begin{proof}[Proof of Theorem \ref{Thm_optimality}]
Set 
\[u:=f+g\in C^k(M),\]
where $f$ is the subsolution given in Theorem \ref{Thm_quadratic}. By \eqref{relation_ktht}, for every $x\in M$ and $t>0$,
\begin{equation}\label{eq_ttumg}
	T_t^-u(x)+c_Ht=\inf_{y\in M}\{u(y)+h_t(y,x)\}+c_Ht=\inf_{y\in M}\{g(y)+h_t^{\hL}(y,x)\}+f(x).
\end{equation}

\noindent\textbf{Step 1.} Let $\ u_\infty:=\lim_{t\to+\infty}(T_t^-u+c_Ht)$. We first show that
\begin{equation}\label{uinf_f}
	  u_\infty(x)=f(x)\quad\text{for all}~x\in \cA_0.
\end{equation}
Indeed, by Theorem \ref{uinftyexpression} and \eqref{relation_ktht}, we have
\begin{equation}\label{q0yf}
	u_\infty(x)=\inf_{y\in M}\{u(y)+h(y,x)\}
	=\inf_{y\in M}\{g(y)+h^{\hL}(y,x)\}+f(x).		
\end{equation}
Since $h^{\hL}\ge 0$ and $g\ge 0$, we obtain $u_\infty\ge f$ on $M$. Next, by Lemma \ref{lem_newL_properties}(ii), 
\begin{equation}\label{eq_zeroLbgam}
\hL(\bgamma(s),\dot{\bgamma}(s))=0,\quad\text{for all}~s.
\end{equation}
It follows that $h^{\hL}(\bgamma(0),\bgamma(s))=0$ for all $s\in[0,T]$. Hence, by \eqref{q0yf}, for each $s\in [0,T]$,
\[
\begin{aligned}
u_\infty(\bgamma(s))
\leq g(\bgamma(0))+h^{\hL}(\bgamma(0),\bgamma(s))+f(\bgamma(s))=f(\bgamma(s)),	
\end{aligned}	
\]
where $ g(\bgamma(0))=0$ follows from \eqref{g_choice}. Therefore, $u_\infty=f$ on $\cA_0$. This proves \eqref{uinf_f}.

\smallskip
\noindent\textbf{Step 2.}
Fix $t\ge T$. We set
\[
\bar{x}_t:=\bgamma(t+T/2)\in \cA_0\qquad\text{and}
\qquad  (N-1)T<t\leq NT\,,
\]
for some integer $N\in \Z_+$. Let $\eta:[0,t]\to M$ be any continuous piecewise $C^1$ curve with $\eta(t)=\bar{x}_t$. Extend it to a curve $\tilde\eta:[0,NT]\to M$ by setting
\[
\tilde\eta(s)=
\begin{cases}
\eta(s),&\quad s\in[0, t],\\
\bgamma(s+T/2),&\quad s\in[t, NT].
\end{cases}
\]
By \eqref{eq_zeroLbgam}, the added segment has zero $\hL$-action. Hence,
\begin{equation*}
\begin{aligned}
\int_0^t\hL(\eta(s),\dot\eta(s))\,ds
=\int_0^{NT}\hL(\tilde\eta(s),\dot{\tilde\eta}(s))\,ds
&\geq \sum_{j=1}^N h_T^{\hL}\Big(\tilde\eta((j-1)T), \tilde\eta(jT)\Big).
\end{aligned}
\end{equation*}
To estimate the right-hand side, set $a_j:=g(\tilde\eta(jT))$ for $j=0,\ldots,N$. By periodicity and \eqref{g_choice}, $a_N=g(\bgamma(T/2))=1$. Moreover, $a_0=g(\eta(0))\in[0,1]$, so $a_0\ge a_0^2$. Combining the preceding estimate with Lemma \ref{lem_one_period_lowerbound} and the Cauchy--Schwarz inequality, we obtain
\begin{align*}
g(\eta(0))+
\int_0^t\hL(\eta(s),\dot\eta(s))\,ds
&\ge a_0^2+\frac1K\sum_{j=1}^N (a_j-a_{j-1})^2\\
&\ge
\frac{\left(a_0+\sum_{j=1}^N (a_j-a_{j-1})\right)^2}{1+KN}=\frac1{1+KN}.
\end{align*}
Taking the infimum over all such curves and using $N<t/T+1$ gives
\begin{equation*}
	\inf_{y\in M} \{g(y)+h_t^\hL(y,\bar{x}_t)\}\geq \frac1{1+KN} >\frac{c}{t},
\qquad\text{for all}~t\geq T,
\end{equation*}
where $c=T/(1+2K)$. Since $\bar{x}_t\in\cA_0$, $u_\infty(\bar{x}_t)=f(\bar{x}_t)$ by \eqref{uinf_f}. This, together with \eqref{eq_ttumg}, yields
\[
\begin{aligned}
\|T_t^-u+c_Ht-u_\infty\|_{C(M)}
\ge T_t^-u(\bar{x}_t)+c_Ht-u_\infty(\bar{x}_t)
&=T_t^-u(\bar{x}_t)+c_Ht-f(\bar{x}_t)\\
&=\inf_{y\in M} \{g(y)+h_t^\hL(y,\bar{x}_t)\}\\
&\ge c\,t^{-1},\qquad t\ge T.
\end{aligned}
\]
Combining this lower bound with Theorem \ref{Thm_semigroup} yields Theorem \ref{Thm_optimality}.
\end{proof}

We conclude this section by proving the lower-bound estimate \eqref{h_lowerbd} in Theorem \ref{Thm_kernel}.

Indeed, let $u$ be the initial datum obtained in Theorem \ref{Thm_optimality}. Combining the lower bound obtained there with \eqref{ineq_Tleqh}, we get
\[
\|h_t+c_Ht-h\|_{C(M\times M)}\geq \|T_t^-u+c_Ht-u_\infty\|_{C(M)}
\geq c\,t^{-1},
\qquad t\geq t_*.
\]
Together with Theorem \ref{thm_equivstatement}, this finally completes the proof of Theorem \ref{Thm_kernel}.

%%%%%%%%%%%%%%%%%%%%%%%%%%%%%%%%%%%%%%%%%%%%%%%%%%%%%%%%%%%%%%%%%%
\vspace{2mm}

\noindent{\bf Acknowledgments.}
~The first author was supported by the National Natural Science Foundation of China (Grant Nos. 12301228 and 12571207) and the Fundamental Research Funds for the Central Universities (Grant No. 020314380038).

\end{document}